\documentclass[12pt,leqno]{article}
\usepackage{amssymb}
\usepackage[mathscr]{eucal}
\usepackage{amsmath,amssymb,latexsym,theorem,bbm}
\usepackage{color}
\usepackage{appendix, url}
\usepackage{graphicx}
\usepackage{subcaption}

\newcommand{\CC}{\mathbb{C}}

\newcommand{\NN}{\mathbb{N}}

\newcommand{\RR}{\mathbb{R}}

\newcommand{\ZZ}{\mathbb{Z}}

\newcommand{\bA}{{\boldsymbol{A}}}

\newcommand{\bB}{{\boldsymbol{B}}}

\newcommand{\be}{{\boldsymbol{e}}}

\newcommand{\bI}{{\boldsymbol{I}}}

\newcommand{\bM}{{\boldsymbol{M}}}

\newcommand{\bN}{{\boldsymbol{N}}}

\newcommand{\bP}{{\boldsymbol{P}}}

\newcommand{\bQ}{{\boldsymbol{Q}}}

\newcommand{\bV}{{\boldsymbol{V}}}

\newcommand{\bx}{{\boldsymbol{x}}}
\newcommand{\bX}{{\boldsymbol{X}}}
\newcommand{\by}{{\boldsymbol{y}}}
\newcommand{\bY}{{\boldsymbol{Y}}}

\newcommand{\bZ}{{\boldsymbol{Z}}}
\newcommand{\bU}{{\boldsymbol{U}}}

\newcommand{\Beta}{{\boldsymbol{\eta}}}

\newcommand{\btheta}{{\boldsymbol{\theta}}}
\newcommand{\bvartheta}{{\boldsymbol{\vartheta}}}

\newcommand{\bzeta}{{\boldsymbol{\zeta}}}
\newcommand{\bzero}{{\boldsymbol{0}}}

\newcommand{\cB}{{\mathcal B}}

\newcommand{\cD}{{\mathcal D}}

\newcommand{\cG}{{\mathcal G}}
\newcommand{\cF}{{\mathcal F}}

\newcommand{\cN}{{\mathcal N}}

\newcommand{\cO}{{\mathcal O}}

\newcommand{\dd}{\mathrm{d}}
\newcommand{\ee}{\mathrm{e}}
\newcommand{\ii}{\mathrm{i}}

\newcommand{\EE}{\operatorname{\mathbb{E}}}

\newcommand{\PP}{{\operatorname{\mathbb{P}}}}

\newcommand{\sign}{\operatorname{sign}}

\newcommand{\hvartheta}{\widehat{\vartheta}}

\newcommand{\vare}{\varepsilon}

\renewcommand{\mid}{\,|\,}

\renewcommand{\leq}{\leqslant}
\renewcommand{\geq}{\geqslant}

\newcommand{\stoch}{\stackrel{\PP}{\longrightarrow}}
\newcommand{\stochG}{\stackrel{\PP_G}{\longrightarrow}}

\newcommand{\stochGeta}{\stackrel{\PP_{G\cap\{\exists\Beta^{-1}\}}}{\longrightarrow}}

\newcommand{\distrP}{\stackrel{\cD(\PP)}{\longrightarrow}}

\newcommand{\distre}{\stackrel{\cD}{=}}

\newcommand{\asP}{\stackrel{ \text{$\PP$-{\rm a.s.}} }{\longrightarrow}}

\newcommand{\aseP}{\stackrel{{\text{$\PP$-a.s.}}}{=}}

\newcommand{\bbone}{\mathbbm{1}}

\newcommand{\proofend}{\hfill\mbox{$\Box$}}

\numberwithin{equation}{section}

\theoremstyle{change} \theorembodyfont{\em}
\newtheorem{Lem}{Lemma.}[section]
\newtheorem{Thm}[Lem]{Theorem.}

\newtheorem{Cor}[Lem]{Corollary.}
\newtheorem{Def}[Lem]{Definition.}

\theorembodyfont{\rm}
\newtheorem{Rem}[Lem]{Remark.}

\begin{document}

\begin{center}
 {\bfseries\Large
   Mixing convergence of LSE for supercritical AR(2)\\[-1mm] processes with Gaussian innovations\\[1mm] using random scaling}

\vspace*{3mm}

 {\sc\large
  M\'aty\'as $\text{Barczy}^{*,\diamond}$,
   Fanni $\text{Ned\'enyi}^{*}$,
  \ \framebox[1.1\width]{Gyula $\text{Pap}$}}

\end{center}

\vskip0.2cm

\noindent
 * HUN-REN–SZTE Analysis and Applications Research Group,
   Bolyai Institute, University of Szeged,
   Aradi v\'ertan\'uk tere 1, H--6720 Szeged, Hungary.

\noindent e--mails: barczy@math.u-szeged.hu (M. Barczy),
                    nfanni@math.u-szeged.hu (F. Ned\'enyi).

\noindent $\diamond$ Corresponding author.

\vskip0.5cm


{\renewcommand{\thefootnote}{}
\footnote{M\'aty\'as Barczy and Fanni Ned\'enyi were supported by the project TKP2021-NVA-09.
Project no.\ TKP2021-NVA-09 has been implemented with the support
 provided by the Ministry of Innovation and Technology of Hungary from the National Research, Development and Innovation Fund,
 financed under the TKP2021-NVA funding scheme.}}

\vspace*{-10mm}

\begin{abstract}
We prove mixing convergence of the least squares estimator of autoregressive parameters for supercritical autoregressive processes of order 2
 with Gaussian innovations having real characteristic roots with different absolute values.
We use an appropriate random scaling such that the limit distribution is a two-dimensional normal distribution
 concentrated on a one-dimensional ray determined by the characteristic root having the larger absolute value.
\end{abstract}

\section{Introduction}
\label{section_intro}

Studying asymptotic behaviour of the Least Squares Estimator (LSE) of AutoRegressive (AR) parameters of AR processes has a long history,
 it goes back at least to Mann and Wald \cite{ManWal}.
Most of the authors proved convergence in distribution of appropriately normalized versions of the LSE in question,
 but one can rarely find other types of convergence in the corresponding limit theorems.
For some AR processes of order \ $p$ \ (AR($p$) processes), Jeganathan \cite[Theorems 9, 14 and 17]{Jeg} proved so-called strong convergence
 (see Jeganathan \cite[Definition 3]{Jeg}) of the LSE in question, and for AR(1) processes,
 H\"ausler and Luschgy \cite[Chapter 9]{HauLus} proved so-called stable and mixing convergence (see Appendix \ref{section_stable}) of the LSE.
Below, we will recall the results of H\"ausler and Luschgy \cite{HauLus} on the asymptotic behaviour of the LSE in question for
 supercritical AR(1) processes with Gaussian innovations.
Note that, since stable (mixing) convergence yields convergence in distribution, the results of H\"ausler and Luschgy \cite[Chapter 9]{HauLus}
 immediately imply convergence in distribution of the LSE in question.
According to our knowledge, results on stable (mixing) convergence of the LSE of AR parameters of higher order AR processes are not available
 in the literature.
In the present paper we consider a supercritical AR process of order 2 with Gaussian innovations having real characteristic roots with different absolute values,
 and we prove mixing convergence of the LSE of its AR parameters using some random normalization.

Limit theorems stating stable (mixing) convergence instead of convergence in distribution are important not only from
 theoretical point of view.
These types of limit theorems have such statistical applications where limit theorems stating only
 convergence in distribution cannot be directly used.
Such a nice application for the description of the asymptotic behaviour of the so-called Harris estimator of the offspring mean of
 a supercritical Galton--Watson process is explained on pages 3 and 4 in H\"ausler and Luschgy \cite{HauLus}.
In the heart of this application there is a generalization of Slutsky's lemma in a way that convergence in distribution is replaced by stable convergence
 and in return the limit of the stochastically convergent sequence can be also random not only a deterministic constant.

Let \ $\ZZ_+$, \ $\NN$, \ $\RR$, \ $\RR_+$ \ and \ $\RR_{++}$ \ denote the set of non-negative integers, positive integers, real numbers,
 non-negative real numbers and positive real numbers, respectively.
Let \ $(Z_n)_{n\in\NN}$ \ be a sequence of independent and identically distributed (i.i.d.) random variables such that
 \ $Z_1$ \ is normally distributed with mean \ $0$ \  and with variance \ $\sigma^2$, \ where \ $\sigma \in \RR_{++}$.

Now we recall the results of H\"ausler and Luschgy \cite[Example 8.10]{HauLus} on mixing convergence of the LSE
 of the AR parameter of a supercritical AR(1) process with Gaussian innovations.
Let \ $Y_0$ \ be a random variable independent of \ $(Z_n)_{n\in\NN}$ \ and suppose that \ $Y_0$ \ has a finite second moment.
Let us consider an AR(1) process \ $(Y_n)_{n\in\ZZ_+}$ \ with Gaussian innovations
 defined by \ $Y_n = \vartheta Y_{n-1} + Z_n$, \ $n\in\NN$, \ where \ $\vartheta\in\RR$.
\ The AR(1) process \ $(Y_n)_{n\in\ZZ_+}$ \ is called subcritical, critical and supercritical if
 \ $\vert \vartheta\vert<1$, \ $\vert \vartheta\vert=1$ \ and \ $\vert \vartheta\vert>1$, \ respectively.
In case of \ $\vert \vartheta\vert>1$, \ i.e., for supercritical AR(1) processes with Gaussian innovations, H\"ausler and Luschgy \cite[Example 8.10]{HauLus}
 proved mixing convergence of \ $\left(\sum_{j=1}^n Y_{j-1}^2\right)^{\frac{1}{2}}(\hvartheta^{(n)} - \vartheta)$ \
 towards a Gaussian distribution with mean \ $0$ \ and variance \ $\sigma^2$ \ as \ $n\to\infty$, \ where
 \[
   \hvartheta^{(n)} = \frac{\sum_{j=1}^n Y_jY_{j-1}}{\sum_{j=1}^n Y_{j-1}^2},
                          \qquad n\in\NN,
 \]
 denotes the LSE of \ $\vartheta$ \ based on the observations \ $Y_0,Y_1,\ldots,Y_n$ \ (provided that the denominator above is not zero).
For results on stable (mixing) convergence of \ $\hvartheta^{(n)}$ \ in case of subcritical and critical AR(1) processes, see
H\"ausler and Luschgy \cite[Theorems 9.1 and 9.3]{HauLus}.

In what follows, let \ $(X_0, X_{-1})^\top$ \ be a random vector with values in \ $\RR^2$ \ independent of \ $(Z_n)_{n\in\NN}$, \
 and suppose that \ $X_0$ \ and \ $X_{-1}$ \ have finite second moments.
Let us consider an AR(2) process \ $(X_n)_{n\geq-1}$ \ with Gaussian innovations defined by
 \begin{equation}\label{AR(2)}
  X_n = \vartheta_1 X_{n-1} + \vartheta_2 X_{n-2} + Z_n , \qquad n \in \NN ,
 \end{equation}
 where \ $(\vartheta_1, \vartheta_2)^\top \in \RR^2$.
\ Note that \eqref{AR(2)} implies
 \begin{equation*}
  \begin{bmatrix} X_n \\ X_{n-1} \end{bmatrix}
  = \bvartheta
    \begin{bmatrix} X_{n-1} \\ X_{n-2} \end{bmatrix}
    + \begin{bmatrix} Z_n \\ 0 \end{bmatrix} , \qquad n \in \NN ,
 \end{equation*}
 and hence
 \begin{equation}\label{X_n_X_n-1}
  \begin{bmatrix} X_n \\ X_{n-1} \end{bmatrix}
  = \bvartheta^k
    \begin{bmatrix} X_{n-k} \\ X_{n-k-1} \end{bmatrix}
    + \sum_{j=n-k+1}^n
       \bvartheta^{n-j}
       \begin{bmatrix} Z_j \\ 0 \end{bmatrix} , \qquad n \in \ZZ_+ , \quad k \in \{0,1, \ldots, n\} ,
 \end{equation}
 where \ $\sum_{j=n+1}^n:=\bzero$ \ and
 \[
   \bvartheta := \begin{bmatrix}
               \vartheta_1 & \vartheta_2 \\
               1 & 0
              \end{bmatrix} .
 \]
By \ $\varrho(\bvartheta)$, \ we denote the spectral radius of \ $\bvartheta$.
\ The AR(2) process \ $(X_n)_{n\geq-1}$ \ given in \eqref{AR(2)} is called subcritical, critical and supercritical
  if \ $\varrho(\bvartheta) < 1$, \ $\varrho(\bvartheta) = 1$ \ and \ $\varrho(\bvartheta) > 1$, \ respectively.
\ We have \ $\varrho(\bvartheta) = \max\{|\lambda_+| , |\lambda_-|\}$, \ where \ $\lambda_+$ \ and \ $\lambda_-$ \  denote the eigenvalues of \ $\bvartheta$ \ given by
 \begin{align}\label{help_lambdak}
   \lambda_+ := \frac{\vartheta_1+\sqrt{\vartheta_1^2+4\vartheta_2}}{2} \in \CC , \qquad
   \lambda_- := \frac{\vartheta_1-\sqrt{\vartheta_1^2+4\vartheta_2}}{2} \in \CC .
 \end{align}
Note that \ $\lambda_+$ \ and \ $\lambda_-$ \ are the (characteristic) roots of the autoregressive (characteristic) polynomial
 \ $x^2 - \vartheta_1 x - \vartheta_2$ \ of the AR(2) process \ $(X_n)_{n\geq-1}$.
\ The process \ $(X_n)_{n\geq-1}$ \ is called explosive if its characteristic polynomial has at least one root outside the unit circle (supercritical case),
 but has no roots on the unit circle, i.e., \ $\varrho(\bvartheta)>1$ \ and \ $\vert \lambda_+\vert\ne1$, \ $\vert \lambda_-\vert\ne1$, \  see Jeganathan \cite[Section 6]{Jeg}.
If both roots \ $\lambda_+$ \ and \ $\lambda_-$ \ are outside the unit circle, i.e., \ $\vert \lambda_+\vert>1$ \ and
 \ $\vert \lambda_-\vert>1$, \ then the process \ $(X_n)_{n\geq-1}$ \ is called purely explosive.
An explosive, but not purely explosive AR(2) process is sometimes called partially explosive.

In this paper we will consider a supercritical AR(2) process with Gaussian innovations
 supposing also that \ $|\lambda_+| \ne |\lambda_-|$.
\ Then \ $\vartheta_1^2 + 4 \vartheta_2 > 0$, \ hence we have \ $\lambda_+ \in \RR$ \ and \ $\lambda_- \in \RR$.
\ Note that this case includes the (purely and partially) explosive case, but it also includes
 for example the case \ $\vert \lambda_+\vert>1$ \ and \ $\lambda_-=1$ \
 (when \ $\vartheta_1>2$, \ $\vartheta_2 = 1 -\vartheta_1$, \ and \ $\lambda_+=\vartheta_1-1$),
 \ i.e., when there is a so-called unit root of the characteristic polynomial of \ $(X_n)_{n\geq-1}$.

We will study the asymptotic behaviour of the LSE of the parameters \ $\vartheta_1$ \ and \ $\vartheta_2$ \  based on the observations \ $X_{-1}$, $X_0$, $X_1$, \ldots, $X_n$ \
 using an appropriate random normalization with the aim to establish mixing convergence (see Appendix \ref{section_stable}) of the LSE in question.
For each \ $n \in \NN$, \ a least squares estimator \ $(\hvartheta_1^{(n)}, \hvartheta_2^{(n)})^\top$ \ of \ $(\vartheta_1, \vartheta_2)^\top$ \ based on the observations \ $X_{-1}$, $X_0$, $X_1$, \ldots, $X_n$ \ can be obtained by minimizing the sum of squares
 \[
   \sum_{k=1}^n (X_k - \vartheta_1 X_{k-1} - \vartheta_2 X_{k-2})^2
 \]
 with respect to \ $(\vartheta_1, \vartheta_2)^\top$ \ over \ $\RR^2$.
\ It is known that for each \ $n \in \NN$ \ with \ $n \geq 3$,
 \ a unique LSE \ $(\hvartheta_1^{(n)}, \hvartheta_2^{(n)})^\top$ \ of \ $(\vartheta_1,\vartheta_2)^\top$ \
 based on the observations \ $X_{-1}, X_0,X_1,\ldots,X_n$ \ exists with probability 1,
 and this LSE has the form given by
 \[
  \begin{bmatrix}
   \hvartheta_1^{(n)} \\
   \hvartheta_2^{(n)}
  \end{bmatrix}
  = \left(\sum_{k=1}^n
     \begin{bmatrix} X_{k-1} \\ X_{k-2} \end{bmatrix}
     \begin{bmatrix} X_{k-1} \\ X_{k-2} \end{bmatrix}^\top\right)^{-1}
    \sum_{k=1}^n
     X_k
     \begin{bmatrix} X_{k-1} \\ X_{k-2} \end{bmatrix} ,
  \qquad n \in \NN,
 \]
 on the event
 \[
    \left\{\det\left(\sum_{k=1}^n
    \begin{bmatrix} X_{k-1} \\ X_{k-2} \end{bmatrix}
    \begin{bmatrix} X_{k-1} \\ X_{k-2} \end{bmatrix}^\top\right) > 0\right\},
 \]
 see, e.g., Lemma \ref{LEMMA_LSE_exist}.

Venkataraman \cite{Ven1}, \cite{Ven2}, \cite{Ven3} and Narasimham \cite{Nar} proved convergence in distribution
 of $(\hvartheta_1^{(n)}, \hvartheta_2^{(n)})^\top$ for explosive AR(2) processes with (not necessarily Gaussian)
 i.i.d.\ innovations using non-random normalizations.
Datta \cite[Theorem 2.2]{Dat} showed convergence in distribution of the LSE of AR parameters of a purely explosive AR process
 of order $p\in\NN$ with (not necessarily Gaussian) i.i.d.\ innovations satisfying a logarithmic moment condition using
 non-random normalization (depending on the unknown parameters to be estimated).
Datta \cite[Corollary 5.1]{Dat} also proved convergence in distribution of the LSE in question for a partially explosive AR process
 of order $p\in\NN$ with (not necessarily Gaussian) i.i.d.\ innovations satisfying a second order moment condition using
 non-random normalization (depending on the unknown parameters).

Touati \cite[Th\'eor\`{e}me 1]{Tou} proved convergence in distribution of the LSE of AR parameters of a purely explosive
 or partially explosive AR process of order \ $p\in\NN$ \ with (not necessarily Gaussian) i.i.d.\ innovations
 having a finite second moment using nonrandom normalization (depending on the unknown parameters) and
 using a kind of random normalization but which also depends on the unknown parameters.
The limit laws were identified as an infinite mixture of random matrices that are not Gaussian in general.
In Section \ref{section_discussion} we give a detailed comparison of our forthcoming results with those of Touati \cite{Tou}.
Recently, Monsour \cite{Mon} has provided a unified approach for proving convergence in distribution of LSE of the AR parameters
 of an AR(p) process of order \ $p\in\NN$ \  with Gaussian innovations.
Part (b) of Theorem 3 in Monsour \cite{Mon} gives a general result on the asymptotic behaviour of the LSE of AR parameters of
 a supercritical AR process of order \ $p\in\NN$ \ with Gaussian innovations
 using a non-random normalization establishing weak convergence of the LSE in question.
Part (c) of Theorem 3 in Monsour \cite{Mon} describes the asymptotic behaviour of the LSE of AR parameters of a general AR process of order
 \ $p\in\NN$ \ with Gaussian innovations using a random normalization proving weak convergence of the LSE in question.
In Section \ref{section_discussion} we give a detailed comparison of our forthcoming results with those of Monsour \cite{Mon}.

Recently, Aknouche \cite{Akn} has studied the asymptotic behaviour of the LSE for some explosive strong periodic AR processes of order $p\in\NN$
 and in case of independent and periodically distributed Gaussian innovations with zero mean, it has been shown that the LSE using an appropriate random scaling
 converges in distribution towards a $p$-dimensional standard normal distribution.

Results on asymptotic behaviour of the LSE of AR parameters of AR processes of order \ $p\in\NN$ \ with not necessarily
  independent innovations are also available in the literature.
Chan and Wei \cite{ChaWei} proved convergence in distribution of the LSE in question for subcritical and critical AR processes of
 order \ $p\in\NN$ \ with innovations that are martingale differences with respect to an increasing sequence of \ $\sigma$-algebras.
Boutahar \cite{Bou} considered AR processes of order \ $p\in\NN$ \ with innovations that form a stationary Gaussian process with
 a regularly varying spectral density, and convergence in distribution of the LSE of AR parameters was established in the subcritical, critical
 and supercritical cases.

In the present paper we consider a supercritical AR(2) process $(X_k)_{k\geq -1}$ given in \eqref{AR(2)}
 with real characteristic roots $\lambda_+,\lambda_-\in\RR$ satisfying $\vert \lambda_+\vert\ne \vert \lambda_-\vert$ and
 we show that
 \begin{align}\label{help_scaled_error}
   &\begin{bmatrix}
     \left(\sum_{k=1}^n X_{k-1}^2\right)^{1/2}
      & \frac{\sum_{k=1}^n X_{k-1} X_{k-2}}{\left(\sum_{k=1}^n X_{k-1}^2\right)^{1/2}} \\[3mm]
     \frac{\sum_{k=1}^n X_{k-1} X_{k-2}}{\left(\sum_{k=1}^n X_{k-2}^2\right)^{1/2}}
      & \left(\sum_{k=1}^n X_{k-2}^2\right)^{1/2}
    \end{bmatrix}
    \begin{bmatrix}
     \hvartheta_1^{(n)} - \vartheta_1 \\
     \hvartheta_2^{(n)} - \vartheta_2
    \end{bmatrix}
 \end{align}
 converges mixing to $\sigma N [ 1 \; \sign(\lambda_1) ]^\top$ as $n\to\infty$, where $\lambda_1$ denotes the
 characteristic root having the larger absolute value, $N$ is a one-dimensional standard normally distributed random variable,
 see Theorem \ref{|lambda_1|>|lambda_2|}.
The limit distribution is in fact a two-dimensional normal distribution concentrated on a one-dimensional ray determined
 by the characteristic root having the larger absolute value.
The random normalization above for $[\hvartheta_1^{(n)} - \vartheta_1,\hvartheta_2^{(n)} - \vartheta_2]^\top$
 seems to be new in the literature, e.g., compared to Touati \cite[Th\'eor\`{e}mes 1 et 2]{Tou} and Monsour \cite[Theorem 3]{Mon}.
Our proof is based on a multidimensional stable limit theorem (see Theorem \ref{MSLTES}, proved in Barczy and Pap \cite{BarPap1})
 which is a multidimensional analogue of the corresponding one-dimensional result in H\"ausler and Luschgy \cite[Theorem 8.2]{HauLus}.
Our proof technique is motivated by that of Theorem 9.2 in H\"ausler and Luschgy \cite{HauLus}, and it is completely different
 from that of Monsour \cite{Mon}.
We note that we tried to prove mixing (or stable-) convergence of $(\hvartheta_1^{(n)},\hvartheta_2^{(n)})$ as $n\to\infty$
 in the supercritical case using a non-random normalization, but our attempts have not been successful so far (for more details,
 see the end of Section \ref{section_discussion}).

We also call the attention to the fact that proving mixing convergence for the LSE of AR parameters
 for supercritical AR(2) processes instead for supercritical AR(1) processes is not an obvious and easy generalization.
As we mentioned above, we had to develop a general multidimensional stable limit theorem, see Barczy and Pap \cite[Theorem 1.4]{BarPap1}
 or Theorem \ref{MSLTES}.
Our assumption that the common distribution of the innovations is Gaussian makes the verification of condition (iv)
 of Theorem \ref{MSLTES} relatively easy (see Step 5 in part (iv) of Theorem \ref{|lambda_1|>|lambda_2|}).
However, we emphasize that, for the application of Theorem \ref{MSLTES}, we do not need Gaussianity of the innovations,
 so, later on, one might study non-Gaussian i.i.d.\ innovations as well.

In Section \ref{section_sim} we illustrate our main result Theorem \ref{|lambda_1|>|lambda_2|} using generated sample paths of a supercritical AR(2) process
started from $(X_0,X_{-1})=(0,0)$.
We choose the AR parameters $\vartheta_1$ and $\vartheta_2$ in a way that the characteristic polynomial of the corresponding AR(2) process
 has two different, real roots and the AR(2) process in question is purely explosive (Case 1), partially explosive (Case 2), supercritical with
 a characteristic root $1$ (Case 3), and supercritical with a characteristic root $-1$ (Case 4), respectively.
In Table \ref{Table1} we summarize these four cases, and we give our choices of parameters.
For every particular choice of $(\vartheta_1,\vartheta_2)$ we generate $1000$ independent $102$-length trajectories of the AR(2) process,
 and for each of them we calculate the corresponding LSE $(\hvartheta_1^{(100)}, \hvartheta_2^{(100)})$
 and the scaled error \eqref{help_scaled_error} of the LSE in question.
For every particular choice of $(\vartheta_1,\vartheta_2)$ we plot the density histogram of the scaled errors of LSEs of the parameters.
We calculate the empirical mean of the LSEs $(\hvartheta_1^{(100)},\hvartheta_2^{(100)})$ based on our 1000 replications.
We also calculate the empirical mean, variance, median, skewness, kurtosis and interquartile range of the scaled errors.
We calculate the empirical covariance between the scaled errors of the LSEs $\hvartheta_1^{(100)}$ and $\hvartheta_2^{(100)}$ as well.
Using Kolmogorov-Smirnov, Pearson's chi-squared, Anderson-Darling and Jarque-Bera tests, we test whether the scaled errors follow a normal distribution.
Our simulation results are in accordance with our theoretical results in Theorem \ref{|lambda_1|>|lambda_2|},
 for a detailed description, see Section \ref{section_sim}.

The paper is structured as follows.
In Section \ref{section_LSE}, we recall the derivation of the LSE $(\hvartheta_1^{(n)}, \hvartheta_2^{(n)})^\top$ of $(\vartheta_1,\vartheta_2)^\top$,
 and a useful decomposition of $(\hvartheta_1^{(n)} - \vartheta_1, \hvartheta_2^{(n)} - \vartheta_2)^\top$ as well, which is used
 in the proof of Theorem \ref{|lambda_1|>|lambda_2|}.
Section \ref{section_AR2_random} contains the precise formulation of our main result Theorem \ref{|lambda_1|>|lambda_2|} together
 with its two corollaries.
In Remark \ref{Rem_confidence} we discuss how one can use Theorem \ref{|lambda_1|>|lambda_2|} for constructing asymptotic
 two-dimensional confidence regions for the AR parameter \ $[\vartheta_1,\vartheta_2]^\top$.
\ Section \ref{section_discussion} is devoted to give a detailed comparison of Theorem \ref{|lambda_1|>|lambda_2|} with the existing results in the literature,
 especially with those of Monsour \cite{Mon} and Touati \cite{Tou}.
Section \ref{section_sim} contains a summary of our simulation results, which confirm our theoretical results.
All the proofs for Sections \ref{section_LSE} and \ref{section_AR2_random} can be found in Section \ref{section_proofs}.
We close the paper with three appendices: we recall the notions of stable and mixing convergence (see Appendix \ref{section_stable}),
 a multidimensional analogue of a one-dimensional stable limit theorem in H\"ausler and Luschgy \cite[Theorem 8.2]{HauLus}
 which was proved in Barczy and Pap \cite[Theorem 1.4]{BarPap1} (see Appendix \ref{section_MSLTES}), and the Lenglart's inequality (see Appendix \ref{section_Lenglart}).

In what follows, we collect the notations used in the paper and not defined so far.
Let $\log^+(x):= \log(x)\bbone_{\{x\geq 1\}} + 0\cdot \bbone_{\{ 0\leq x < 1\}}$ for $x\in\RR_+$.
The Borel $\sigma$-algebra on $\RR^d$ is denoted by $\cB(\RR^d)$, where $d\in\NN$.
Almost sure convergence with respect to a probability measure $\PP$,
 convergence in a probability measure $\PP$ and in distribution under a probability measure $\PP$ will be denoted by
 $\asP$, \ $\stoch$ and $\distrP$, respectively.
For an event $A$ with $\PP(A) > 0$,
 let $\PP_A(\cdot) := \PP(\cdot\mid A) = \PP(\cdot \cap A) / \PP(A)$ denote
 the conditional probability measure given $A$.
Let $\EE_\PP$ denote the expectation under the probability measure $\PP$.
Almost sure equality under a probability measure $\PP$ and equality in distribution will be denoted by $\aseP$
 and $\distre$, respectively.
Every random variable will be defined on a suitable probability space $(\Omega,\cF,\PP)$.
For a random variable $\xi:\Omega\to\RR^d$, its distribution on $(\RR^d,\cB(\RR^d))$ under $\PP$ is denoted by $\PP^\xi$.
By $\|\bx\|$ and $\|\bA\|$, we denote the Euclidean norm of a vector $\bx \in \RR^d$ and the induced matrix norm of a matrix
 $\bA \in \RR^{d\times d}$, respectively.
By $\langle \bx,\by\rangle$, we denote the Euclidean inner product of vectors $\bx,\by\in\RR^d$.
The null vector and the null matrix will be denoted by $\bzero$.
Moreover, $\bI_d \in \RR^{d\times d}$ denotes the identity matrix, and $\be_1$, \ldots, $\be_d$ denotes the natural bases in $\RR^d$.
For a symmetric and positive semidefinite matrix $\bA \in \RR^{d\times d}$,  its unique symmetric, positive semidefinite square root
 is denoted by \ $\bA^{1/2}$.
If $\bV \in \RR^{d\times d}$ is symmetric and positive semidefinite, then
 $\cN_d(\bzero, \bV)$ denotes the $d$-dimensional normal distribution with
 mean vector $\bzero\in\RR^d$ and covariance matrix $\bV$.
In case of $d=1$, instead of $\cN_1$ we simply write $\cN$.

\section{Preliminaries on LSE}
\label{section_LSE}

First, we define the notion of a LSE of $(\vartheta_1, \vartheta_2)^\top$ based on discrete time observations, and then we derive an explicit expression for it.
For each $n \in \NN$, a LSE $(\hvartheta_1^{(n)}, \hvartheta_2^{(n)})^\top$ of
 $(\vartheta_1, \vartheta_2)^\top$ based on the observations $X_{-1}, X_0, X_1, \ldots, X_n$ can be obtained by minimizing the sum of squares
 \[
   \sum_{k=1}^n (X_k - \vartheta_1 X_{k-1} - \vartheta_2 X_{k-2})^2
 \]
 with respect to $(\vartheta_1, \vartheta_2)^\top$ over $\RR^2$.
For each $n \in \NN$, we define the
 function $Q_n : \RR^{n+2} \times \RR^2 \to \RR$ by
 \[
   Q_n(x_{-1}, x_0, x_1, \ldots, x_n ; \vartheta_1, \vartheta_2)
   := \sum_{k=1}^n (x_k - \vartheta_1 x_{k-1} - \vartheta_2 x_{k-2})^2
 \]
 for all \ $(x_{-1}, x_0, x_1, \ldots, x_n)^\top \in \RR^{n+2}$ \ and \ $(\vartheta_1, \vartheta_2)^\top \in \RR^2$.
\ By definition, for each \ $n \in \NN$, \ a least squares estimator of \ $(\vartheta_1, \vartheta_2)^\top$ \ is
 a measurable function \ $F_n : \RR^{n+2} \to \RR^2$ \ such that
 \begin{align*}
  Q_n(x_{-1}, x_0, x_1, \ldots, x_n; F_n(x_{-1}, x_0, x_1, \ldots, x_n))
   = \inf_{(\vartheta_1, \vartheta_2)^\top \in \RR^2} Q_n(x_{-1},x_0, x_1, \ldots, x_n ; \vartheta_1, \vartheta_2)
 \end{align*}
 for all \ $(x_{-1}, x_0, x_1, \ldots, x_n)^\top \in \RR^{n+2}$.
\ Next, we give the solutions of this extremum problem.

\begin{Lem}\label{LSE1}
For each \ $n \in \NN$, \ any least squares estimator of \ $(\vartheta_1, \vartheta_2)^\top$ \ is a measurable function \ $F_n : \RR^{n+2} \to \RR^2$ \ for which
 \begin{equation}\label{LSE_vartheta_1_vartheta_2}
  F_n(x_{-1}, x_0, x_1, \ldots, x_n)
  = G_n(x_{-1}, x_0, x_1, \ldots, x_n)^{-1} H_n(x_{-1}, x_0, x_1, \ldots, x_n)
 \end{equation}
 on the set
 \[
   D_n := \bigl\{ (x_{-1}, x_0, x_1, \ldots, x_n)^\top \in \RR^{n+2}
                  : \det(G_n(x_{-1}, x_0, x_1, \ldots, x_n)) > 0 \bigr\} ,
 \]
 where
 \[
   G_n(x_{-1}, x_0, x_1, \ldots, x_n)
   := \sum_{k=1}^n
       \begin{bmatrix} x_{k-1} \\ x_{k-2} \end{bmatrix}
       \begin{bmatrix} x_{k-1} \\ x_{k-2} \end{bmatrix}^\top , \qquad
   H_n(x_{-1}, x_0, x_1, \ldots, x_n)
   := \sum_{k=1}^n
       x_k
       \begin{bmatrix} x_{k-1} \\ x_{k-2} \end{bmatrix} .
 \]
\end{Lem}

The next result is about the unique existence of \ $(\hvartheta_1^{(n)}, \hvartheta_2^{(n)})^\top$.

\begin{Lem}\label{LEMMA_LSE_exist}
Let \ $(X_k)_{k\geq-1}$ \ be an AR(2) process given in \eqref{AR(2)} such that \ $(Z_n)_{n\in\NN}$ \ is a sequence of i.i.d.\ random variables with \ $Z_1 \distre \cN(0, \sigma^2)$, \ where \ $\sigma \in \RR_{++}$, \ and \ $(X_0, X_{-1})^\top$ \ is a random vector with values in \ $\RR^2$ \ independent of \ $(Z_n)_{n\in\NN}$ \ and with \ $\EE_\PP(X_0^2) < \infty$, \ $\EE_\PP(X_{-1}^2) < \infty$.
\ Then for each \ $n \in \NN$ \ with \ $n \geq 3$, \ we have \ $\PP(\Omega_n) = 1$ \ for the event \ $\Omega_n$ \ given by
 \begin{equation}\label{Omega_n}
  \Omega_n
  := \left\{\det\left(\sum_{k=1}^n
     \begin{bmatrix} X_{k-1} \\ X_{k-2} \end{bmatrix}
     \begin{bmatrix} X_{k-1} \\ X_{k-2} \end{bmatrix}^\top\right) > 0\right\} ,
 \end{equation}
 and hence a unique least squares estimator \ $(\hvartheta_1^{(n)}, \hvartheta_2^{(n)})^\top$ \ of \ $(\vartheta_1,\vartheta_2)^\top$ \
 based on the observations \ $X_{-1}, X_0,X_1,\ldots,X_n$ \ exists with probability 1,
 and this least squares estimator has the form given by
 \begin{equation}\label{LSE}
  \begin{bmatrix}
   \hvartheta_1^{(n)} \\
   \hvartheta_2^{(n)}
  \end{bmatrix}
  = \left(\sum_{k=1}^n
     \begin{bmatrix} X_{k-1} \\ X_{k-2} \end{bmatrix}
     \begin{bmatrix} X_{k-1} \\ X_{k-2} \end{bmatrix}^\top\right)^{-1}
    \sum_{k=1}^n
     X_k
     \begin{bmatrix} X_{k-1} \\ X_{k-2} \end{bmatrix} ,
  \qquad n \in \NN,
 \end{equation}
 on the event \ $\Omega_n$.
\end{Lem}

By Lemma \ref{LEMMA_LSE_exist}, for each \ $n \in \NN$ \ with \ $n\geq 3$, \ on the event \ $\Omega_n$ \ having probability 1, we have
 \begin{align*}
  \begin{bmatrix}
   \hvartheta_1^{(n)} - \vartheta_1 \\
   \hvartheta_2^{(n)} - \vartheta_2
  \end{bmatrix}
  &= \left(\sum_{k=1}^n
      \begin{bmatrix} X_{k-1} \\ X_{k-2} \end{bmatrix}
      \begin{bmatrix} X_{k-1} \\ X_{k-2} \end{bmatrix}^\top\right)^{-1}
     \sum_{k=1}^n
      X_k
      \begin{bmatrix} X_{k-1} \\ X_{k-2} \end{bmatrix}
     - \begin{bmatrix} \vartheta_1 \\ \vartheta_2 \end{bmatrix} \\
  &= \left(\sum_{k=1}^n
      \begin{bmatrix} X_{k-1} \\ X_{k-2} \end{bmatrix}
      \begin{bmatrix} X_{k-1} \\ X_{k-2} \end{bmatrix}^\top\right)^{-1}
     \sum_{k=1}^n
      X_k
      \begin{bmatrix} X_{k-1} \\ X_{k-2} \end{bmatrix} \\
  &\quad
     - \left(\sum_{k=1}^n
        \begin{bmatrix} X_{k-1} \\ X_{k-2} \end{bmatrix}
        \begin{bmatrix} X_{k-1} \\ X_{k-2} \end{bmatrix}^\top\right)^{-1}
       \left(\sum_{k=1}^n
        \begin{bmatrix} X_{k-1} \\ X_{k-2} \end{bmatrix}
        \begin{bmatrix} X_{k-1} \\ X_{k-2} \end{bmatrix}^\top\right)
       \begin{bmatrix} \vartheta_1 \\ \vartheta_2 \end{bmatrix}\\
  &= \left(\sum_{k=1}^n
      \begin{bmatrix} X_{k-1} \\ X_{k-2} \end{bmatrix}
      \begin{bmatrix} X_{k-1} \\ X_{k-2} \end{bmatrix}^\top\right)^{-1}
     \sum_{k=1}^n
      (X_k - \vartheta_1 X_{k-1} - \vartheta_2 X_{k-2})
      \begin{bmatrix} X_{k-1} \\ X_{k-2} \end{bmatrix} \\
  &= \left(\sum_{k=1}^n
      \begin{bmatrix} X_{k-1} \\ X_{k-2} \end{bmatrix}
      \begin{bmatrix} X_{k-1} \\ X_{k-2} \end{bmatrix}^\top\right)^{-1}
     \sum_{k=1}^n
      Z_k
      \begin{bmatrix} X_{k-1} \\ X_{k-2} \end{bmatrix},
 \end{align*}
 where the last equality follows by \eqref{AR(2)}.
Hence
 \begin{align}\label{LSE_diff_formula}
   \begin{bmatrix}
   \hvartheta_1^{(n)} - \vartheta_1 \\
   \hvartheta_2^{(n)} - \vartheta_2
   \end{bmatrix}
   = \langle \bM\rangle_n^{-1} \bM_n,
  \end{align}
 where
 \begin{align}\label{Def_Mn}
   \bM_n := \sigma^{-2} \sum_{k=1}^n Z_k \begin{bmatrix} X_{k-1} \\ X_{k-2} \end{bmatrix} , \qquad n \in \NN,
            \qquad \text{with \ $\bM_0 := \bzero$,}
 \end{align}
 is a square integrable martingale with respect to the filtration \ $(\cF_n)_{n\in\ZZ_+}$, \
 where \ $\cF_n := \sigma(X_{-1}, X_0, X_1,\ldots, X_n) = \sigma(X_{-1}, X_0, Z_1,\ldots, Z_n)$, \ $n\in\ZZ_+$, \
 and \ $(\langle \bM\rangle_n)_{n\in\ZZ_+}$ \ is its quadratic characteristic process given by
 \begin{equation}\label{langle_bM_rangle_n}
 \langle \bM\rangle_n
   =  \sigma^{-2}
      \sum_{k=1}^n
       \begin{bmatrix} X_{k-1} \\ X_{k-2} \end{bmatrix}
       \begin{bmatrix} X_{k-1} \\ X_{k-2} \end{bmatrix}^\top , \qquad n\in\NN,
       \qquad \text{with \ $\langle \bM\rangle_0:=\bzero$.}
 \end{equation}
Indeed, using the definition of a quadratic characteristic process (see, e.g., H\"ausler and Luschgy \cite[page 193]{HauLus} or Corollary \ref{Lenglart_Y}),
 for each \ $n \in \NN$, \ we have
 \begin{equation*}
  \begin{aligned}
   \langle \bM\rangle_n
   &:= \sum_{k=1}^n \EE_\PP((\bM_k - \bM_{k-1}) (\bM_k - \bM_{k-1})^\top \mid \cF_{k-1}) \\
   &= \sum_{k=1}^n
       \EE_\PP\left(\biggl(\sigma^{-2} Z_k \begin{bmatrix} X_{k-1} \\ X_{k-2} \end{bmatrix}\biggr)
                \biggl(\sigma^{-2} Z_k \begin{bmatrix} X_{k-1} \\ X_{k-2} \end{bmatrix}\biggr)^\top
                \,\bigg|\, \cF_{k-1}\right) \\
  &= \sigma^{-4}
      \sum_{k=1}^n
       (\EE_\PP(Z_k^2)) \begin{bmatrix} X_{k-1} \\ X_{k-2} \end{bmatrix}
       \begin{bmatrix} X_{k-1} \\ X_{k-2} \end{bmatrix}^\top
    = \sigma^{-2}
      \sum_{k=1}^n
       \begin{bmatrix} X_{k-1} \\ X_{k-2} \end{bmatrix}
       \begin{bmatrix} X_{k-1} \\ X_{k-2} \end{bmatrix}^\top .
  \end{aligned}
 \end{equation*}

As a consequence of \eqref{langle_bM_rangle_n}, we have \ $\Omega_n = \{\det(\langle \bM\rangle_n)>0\}$, \ where \ $\Omega_n$ \ is given in \eqref{Omega_n}.

\section{Mixing convergence of LSE using random scaling}
\label{section_AR2_random}

Using an appropriate random scaling we prove mixing convergence of the LSE \ $(\hvartheta_1^{(n)}, \hvartheta_2^{(n)})^\top$ \
 given in \eqref{LSE} as \ $n\to\infty$ \ in supercritical cases with real characteristic roots having different absolute values.

Considering a supercritical AR(2) process given in \eqref{AR(2)} such that \ $|\lambda_+| \ne |\lambda_-|$, \ let us introduce the notation
 \[
   (\lambda_1, \lambda_2)
   := \begin{cases}
       (\lambda_+, \lambda_-) & \text{if \ $|\lambda_+| > |\lambda_-|$,} \\
       (\lambda_-, \lambda_+) & \text{if \ $|\lambda_-| > |\lambda_+|$,}
      \end{cases}
 \]
 where \ $\lambda_+$ \ and \ $\lambda_-$ \ are given in \eqref{help_lambdak}, i.e., \ $\lambda_1$ \ and \ $\lambda_2$ \ is the characteristic root
 having the larger and smaller absolute value, respectively.

Recall that \ $\cF_n = \sigma(X_{-1},X_0, X_1 , \dots, X_n)$, $n\in\ZZ_+$, \ and let
 \ $\cF_\infty:= \sigma(\bigcup_{n=0}^\infty \cF_n)$.
We suppose that the probability space \ $(\Omega,\cF_\infty,\PP)$ \ is complete,
 i.e., for all \ $S \in \cF_\infty$ \ with \ $\PP(S) = 0$ \ and all subsets \ $U\subset S$,
 \ we have \ $U\in\cF_\infty$.

\begin{Thm}\label{|lambda_1|>|lambda_2|}
Let \ $(X_k)_{k\geq-1}$ \ be an AR(2) process given in \eqref{AR(2)} such that \ $(Z_n)_{n\in\NN}$ \ is a sequence of i.i.d.\
 random variables with \ $Z_1 \distre \cN(0, \sigma^2)$, \ where \ $\sigma \in \RR_{++}$, \ and \ $(X_0, X_{-1})^\top$ \
 is a random vector with values in \ $\RR^2$ \ independent of \ $(Z_n)_{n\in\NN}$ \ and with \ $\EE_\PP(X_0^2) < \infty$, \ $\EE_\PP(X_{-1}^2) < \infty$.
\ Suppose that the autoregressive (characteristic) polynomial \ $x^2 - \vartheta_1 x - \vartheta_2$ \ of
 \ $(X_k)_{k\geq-1}$ \ has real roots \ $\lambda_1$ \ and \ $\lambda_2$ \ with \ $|\lambda_1| > |\lambda_2|$ \ and \ $|\lambda_1| > 1$.
\ Then
 \begin{equation}\label{conv_|lambda_1|>|lambda_2|}
  \begin{aligned}
   &\begin{bmatrix}
     \left(\sum_{k=1}^n X_{k-1}^2\right)^{1/2}
      & \frac{\sum_{k=1}^n X_{k-1} X_{k-2}}{\left(\sum_{k=1}^n X_{k-1}^2\right)^{1/2}} \\[3mm]
     \frac{\sum_{k=1}^n X_{k-1} X_{k-2}}{\left(\sum_{k=1}^n X_{k-2}^2\right)^{1/2}}
      & \left(\sum_{k=1}^n X_{k-2}^2\right)^{1/2}
    \end{bmatrix}
    \begin{bmatrix}
     \hvartheta_1^{(n)} - \vartheta_1 \\
     \hvartheta_2^{(n)} - \vartheta_2
    \end{bmatrix}
   \to \sigma N \begin{bmatrix} 1 \\ \sign(\lambda_1) \end{bmatrix} \qquad
    \text{$\cF_\infty$-mixing}
 \end{aligned}
 \end{equation}
 as \ $n \to \infty$, \ where \ $N$ \ is a standard normally distributed random variable being \ $\PP$-independent of \ $\cF_\infty$.
\end{Thm}

Note that the random scaling matrix in \eqref{conv_|lambda_1|>|lambda_2|} is well-defined with probability 1 for each
 \ $n\in\NN$ \ with \ $n\geq 3$, \ since \ $X_1$ \ is absolutely continuous yielding
 \ $\PP(\sum_{k=1}^n X_{k-1}^2 > 0)=1$ \ and \ $\PP(\sum_{k=1}^n X_{k-2}^2 > 0)=1$.

In the next remark we give another representation of the law of the limit random variable
in \eqref{conv_|lambda_1|>|lambda_2|} in Theorem \ref{|lambda_1|>|lambda_2|}.

\begin{Rem}
By Step 6 of the proof of Theorem \ref{|lambda_1|>|lambda_2|}, the law of the limit random variable in \eqref{conv_|lambda_1|>|lambda_2|} can be
 represented as the law of the \ $\PP$-almost surely convergent series \ $\sigma\sum_{j=0}^\infty \lambda_1^{-j} \bN_j$, \ where
 \ $(\bN_j)_{j\in\ZZ_+}$ \ is a sequence of \ $\PP$-independent and identically distributed \ $\RR^2$-valued random vectors being
 \ $\PP$-independent of \ $\cF_\infty$ \ such that the law of \ $\bN_0$ \ under \ $\PP$ \ is
 \[
        \cN_2\left(\bzero, \frac{\lambda_1^2-1}{\lambda_1^2}
                   \begin{bmatrix}
                    1 & \sign(\lambda_1) \\
                    \sign(\lambda_1) & 1
                   \end{bmatrix}\right).
 \]
\proofend
\end{Rem}

Next, we formulate two corollaries of the proof of Theorem \ref{|lambda_1|>|lambda_2|} about the asymptotic behaviour of
 \ $\bM_n$ \ (given in \eqref{Def_Mn}) and its quadratic characteristic process \ $\langle \bM\rangle_n$ \ (given in \eqref{langle_bM_rangle_n})
  as \ $n\to\infty$, \ proving stable convergence and \ $\PP$-almost sure convergence, respectively.
These two results can be interesting on their own rights as well.

\begin{Cor}\label{Cor1}
Under the conditions of Theorem \ref{|lambda_1|>|lambda_2|},
 for the process \ $(\bM_n)_{n\in\ZZ_+}$ \ defined in \eqref{Def_Mn}, we have
 \[
   \lambda_1^{-n}\bM_n \to \Beta \sum_{j=0}^\infty \lambda_1^{-j}\bN_j
                           \distre \Beta N \begin{bmatrix}
                                             1 \\
                                             \sign(\lambda_1) \\
                                           \end{bmatrix}
    \qquad \text{$\cF_\infty$-stably as \ $n \to \infty$,}
 \]
 where
 \begin{itemize}
   \item the random variable \ $\Beta$ \ is given by
          \[
          \Beta := \frac{|Y|}{\sigma\sqrt{\lambda_1^2-1}}
            \begin{bmatrix}
             1 & 0 \\
             0 & |\lambda_1|^{-1}
            \end{bmatrix}
         \]
         with
         \begin{align}\label{conv_X_n}
             Y := \frac{\lambda_1}{\lambda_1-\lambda_2} (X_0 - \lambda_2 X_{-1})
                  + \frac{\lambda_1}{\lambda_1-\lambda_2} \sum_{j=1}^\infty \lambda_1^{-j} Z_j ,
         \end{align}
         where the series is absolute convergent \ $\PP$-almost surely,
   \item $(\bN_j)_{j\in\ZZ_+}$ \ is a sequence of \ $\PP$-independent and identically distributed \ $\RR^2$-valued random
          vectors being \ $\PP$-independent of \ $\cF_\infty$ \ such that the law of \ $\bN_0$ \ under \ $\PP$ \ is
         \[
              \cN_2\left(\bzero, \frac{\lambda_1^2-1}{\lambda_1^2}
                     \begin{bmatrix}
                      1 & \sign(\lambda_1) \\
                      \sign(\lambda_1) & 1
                    \end{bmatrix}\right),
       \]
   \item the series \ $\sum_{j=0}^\infty \lambda_1^{-j}\bN_j$ \ converges \ $\PP$-almost surely,
   \item $N$ \ is \ $\PP$-independent of \ $\cF_\infty$ \ with \ $N\distre \cN(0,1)$.
         \ Consequently, \ $N$ \ and \ $\Beta$ \ are \ $\PP$-independent, since \ $\Beta$ \ is \ $\cF_\infty$-measurable.
 \end{itemize}
\end{Cor}

Note that the \ $\PP$-almost sure absolute convergence of the series \ $\sum_{j=1}^\infty \lambda_1^{-j} Z_j$ \ in the definition
 of \ $Y$ \ in \eqref{conv_X_n} follows, e.g., from Lemma 8.1 in H\"ausler and Luschgy \cite{HauLus}, since \ $\EE_\PP(\log^+(\vert Z_1\vert))<\infty$ \
 and \ $\vert \lambda_1\vert>1$, \ $\lambda_1\in\RR$.
\ Remark also that \ $Y$ \ is \ $\cF_\infty$-measurable, since the series \ $\sum_{j=1}^\infty \lambda_1^{-j} Z_j$ \
 converges \ $\PP$-a.s., \ $\sum_{j=1}^n \lambda_1^{-j} Z_j$ \ is \ $\cF_\infty$-measurable for all \ $n\in\NN$ \
 (due to \ $\cF_n=\sigma(X_{-1},X_0,Z_1,\ldots,Z_n)$, \ $n\in\NN$), and the probability space \ $(\Omega,\cF_\infty,\PP)$ \ is complete.

Next, we describe the asymptotic behaviour of \ $\langle \bM\rangle_n$ \ (given in \eqref{langle_bM_rangle_n}) as \ $n\to\infty$.

\begin{Cor}\label{Cor_detMn}
Under the conditions of Theorem \ref{|lambda_1|>|lambda_2|}, for the process \ $(\langle \bM\rangle_n)_{n\in\ZZ_+}$ \ given in \eqref{langle_bM_rangle_n},
 we have
 \begin{align}\label{help_deterministic_scale1}
  \begin{split}
  \lambda_1^{-2n}
  \langle \bM\rangle_n
  = \frac{\lambda_1^{-2n}}{\sigma^2}
     \begin{bmatrix}
      \sum_{k=1}^n X_{k-1}^2 & \sum_{k=1}^n X_{k-1} X_{k-2} \\
      \sum_{k=1}^n X_{k-1} X_{k-2} & \sum_{k=1}^n X_{k-2}^2
     \end{bmatrix}
  \asP \frac{Y^2}{(\lambda_1^2-1)\sigma^2}
       \begin{bmatrix}
        1 & \lambda_1^{-1} \\
        \lambda_1^{-1} & \lambda_1^{-2}
       \end{bmatrix}
  \end{split}
 \end{align}
 as \ $n \to \infty$, \ where \ $Y$ \ is given in \eqref{conv_X_n}.
Especially, \ $\lambda_1^{-4n} \det(\langle \bM\rangle_n) \asP 0$ \ as \ $n \to \infty$.
\end{Cor}

In the next remark we discuss how one can use Theorem \ref{|lambda_1|>|lambda_2|} for constructing asymptotic confidence regions for
 the AR parameter \ $[\vartheta_1,\vartheta_2]^\top$.

\begin{Rem}\label{Rem_confidence}
Let us suppose that the conditions of Theorem \ref{|lambda_1|>|lambda_2|} hold.
Then, using \eqref{conv_|lambda_1|>|lambda_2|}, the fact that mixing convergence yields convergence in distribution
 and the continuous mapping theorem, we get
 \begin{align*}
   &\begin{bmatrix}
     \hvartheta_1^{(n)} - \vartheta_1 \\
     \hvartheta_2^{(n)} - \vartheta_2
   \end{bmatrix}^\top
   \begin{bmatrix}
     \left(\sum_{k=1}^n X_{k-1}^2\right)^{1/2}
      & \frac{\sum_{k=1}^n X_{k-1} X_{k-2}}{\left(\sum_{k=1}^n X_{k-1}^2\right)^{1/2}} \\[3mm]
     \frac{\sum_{k=1}^n X_{k-1} X_{k-2}}{\left(\sum_{k=1}^n X_{k-2}^2\right)^{1/2}}
      & \left(\sum_{k=1}^n X_{k-2}^2\right)^{1/2}
    \end{bmatrix}^\top\\
   &\times
   \begin{bmatrix}
     \left(\sum_{k=1}^n X_{k-1}^2\right)^{1/2}
      & \frac{\sum_{k=1}^n X_{k-1} X_{k-2}}{\left(\sum_{k=1}^n X_{k-1}^2\right)^{1/2}} \\[3mm]
     \frac{\sum_{k=1}^n X_{k-1} X_{k-2}}{\left(\sum_{k=1}^n X_{k-2}^2\right)^{1/2}}
      & \left(\sum_{k=1}^n X_{k-2}^2\right)^{1/2}
    \end{bmatrix}
    \begin{bmatrix}
     \hvartheta_1^{(n)} - \vartheta_1 \\
     \hvartheta_2^{(n)} - \vartheta_2
   \end{bmatrix}
  \distrP \sigma^2 N^2
          \begin{bmatrix}
           1 \\
            \sign(\lambda_1) \\
          \end{bmatrix}^\top
          \begin{bmatrix}
           1 \\
            \sign(\lambda_1) \\
          \end{bmatrix}
 \end{align*}
 as \ $n\to\infty$, \ where \ $N\distre \cN(0,1)$, \ yielding that \ $N^2$ \ has a chi-squared distribution with 2 degrees of freedom.
Consequently,
 \begin{align*}
  & \begin{bmatrix}
     \hvartheta_1^{(n)} - \vartheta_1 \\
     \hvartheta_2^{(n)} - \vartheta_2
   \end{bmatrix}^\top
     \begin{bmatrix}
       \sum_{k=1}^n X_{k-1}^2 + \frac{\left(\sum_{k=1}^n X_{k-1}X_{k-2}\right)^2}{\sum_{k=1}^n X_{k-2}^2} &  2 \sum_{k=1}^n X_{k-1}X_{k-2} \\
        2 \sum_{k=1}^n X_{k-1}X_{k-2}  &  \sum_{k=1}^n X_{k-2}^2 + \frac{\left(\sum_{k=1}^n X_{k-1}X_{k-2}\right)^2}{\sum_{k=1}^n X_{k-1}^2} \\
   \end{bmatrix}
  \begin{bmatrix}
     \hvartheta_1^{(n)} - \vartheta_1 \\
     \hvartheta_2^{(n)} - \vartheta_2
   \end{bmatrix} \\
 &\qquad \distrP 2 \sigma^2 N^2 \qquad \text{as \ $n\to\infty$.}
 \end{align*}
Given \ $\alpha\in(0,1)$, \ let \ $\chi_2^2(1-\alpha)$ \ denote the \ $(1-\alpha)$-quantile of the chi-squared distribution
 with 2 degrees of freedom.
Thus, under the conditions of Theorem \ref{|lambda_1|>|lambda_2|}, for large sample size \ $n$, \
 the region
 \begin{align*}
   \Bigg\{
     \begin{bmatrix}
       u \\
       v \\
     \end{bmatrix}\in\RR^2
   &:  \begin{bmatrix}
     \hvartheta_1^{(n)} - u \\
     \hvartheta_2^{(n)} - v
   \end{bmatrix}^\top
     \begin{bmatrix}
       \sum_{k=1}^n X_{k-1}^2 + \frac{\left(\sum_{k=1}^n X_{k-1}X_{k-2}\right)^2}{\sum_{k=1}^n X_{k-2}^2} &  2 \sum_{k=1}^n X_{k-1}X_{k-2} \\
        2 \sum_{k=1}^n X_{k-1}X_{k-2}  &  \sum_{k=1}^n X_{k-2}^2 + \frac{\left(\sum_{k=1}^n X_{k-1}X_{k-2}\right)^2}{\sum_{k=1}^n X_{k-1}^2} \\
   \end{bmatrix}\\
  &\;\;\times \begin{bmatrix}
     \hvartheta_1^{(n)} - u \\
     \hvartheta_2^{(n)} - v
   \end{bmatrix}
  \leq   2 \sigma^2 \chi_2^2(1-\alpha)
  \Bigg\}
 \end{align*}
 contains the (true) AR parameter \ $[\vartheta_1,\vartheta_2]^\top$ \ with probability close to \ $1-\alpha$.
\proofend
\end{Rem}

\section{Discussion on related results}\label{section_discussion}

First, we give a detailed comparison of Theorem \ref{|lambda_1|>|lambda_2|} (of the present paper) and part (c) of Theorem 3 in Monsour \cite{Mon}
 specialized to the two-dimensional case.
We will distinguish two supercritical cases, namely, without a unit root and with a unit root.
Then we compare Theorem \ref{|lambda_1|>|lambda_2|} (of the present paper) and Th\'eor\`{e}me 1 in Touati \cite{Tou} specialized to the two-dimensional case.

First of all, we note that we are not convinced that the proof of part (c) of Theorem 3 in Monsour \cite{Mon} is complete/correct.
Namely, with the notations of Monsour \cite{Mon}, we do not understand how the weak convergence \ $\cO_n^\top \cO_n\distrP \bI_p$ \
 as \ $n\to\infty$ \ yields that \ $\cO_n$ \ converges in distribution as \ $n\to\infty$, \ where
 \ $\cO_n:=(VG\sum_{j=1}^n \widetilde X_{j-1} \widetilde X_{j-1}^\top G^\top V^\top)^{-1/2} VG R_n^{-1} H^{1/2}$, \ $n\in\NN$, \
 and we cannot see why \ $(\cO_n, H^{-1/2}R_n \sum_{j=1}^n \widetilde X_{j-1}\vare_j)$ \ converges in distribution as \ $n\to\infty$ \
 (which is implicitly used in the proof in question).

Let us suppose that the conditions of Theorem \ref{|lambda_1|>|lambda_2|} hold together with \ $\vert \lambda_2\vert\ne 1$, \ i.e.,
 \ $\lambda_1,\lambda_2\in\RR$ \ and either \ $\vert \lambda_1\vert > \vert \lambda_2\vert>1$ \ (purely explosive case) or
 \ $\vert \lambda_1\vert > 1 > \vert \lambda_2\vert$ \ (partially explosive case).
For simplicity, let us suppose that \ $\sigma=1$.
\ Provided that the proof of part (c) of Theorem 3 in Monsour \cite{Mon} is complete, in our considered case it states that
 \begin{align}\label{help_rconv_2dnorm}
  \begin{bmatrix}
     \sum_{k=1}^n X_{k-1}^2 & \sum_{k=1}^n X_{k-1} X_{k-2} \\[1mm]
     \sum_{k=1}^n X_{k-1} X_{k-2} & \sum_{k=1}^n X_{k-2}^2
    \end{bmatrix}^{1/2}
    \begin{bmatrix}
     \hvartheta_1^{(n)} - \vartheta_1 \\
     \hvartheta_2^{(n)} - \vartheta_2
    \end{bmatrix}
   \distrP  \cN_2(\bzero, \bI_2) \qquad \text{as \ $n \to \infty$.}
 \end{align}
We note that in part (c) of Theorem 3 in Monsour \cite{Mon}, the limit law is represented in a more complicated form,
 but from its proof it turns out that in the special case of no characteristic roots on the unit circle
 (as in our considered case) the limit law can be represented in the form \ $\cO\, (\zeta_1,\zeta_2)^\top$,
 \ where \ $\cO$ \ is some random \ $2\times 2$ orthogonal matrix independent of \ $(\zeta_1,\zeta_2)^\top$ \ having a \ $2$-dimensional
 standard normal distribution (see page 305 in Monsour \cite{Mon}).
In this case, the limit law in question is indeed \ $\cN_2(\bzero,\bI_2)$, \ since the law of a (multidimensional) standard normally
 distributed random variable is invariant under an orthogonal transformation, and the random orthogonal matrix \ $\cO$ \ and the
 standard normally distributed (two-dimensional) random vector \ $(\zeta_1,\zeta_2)^\top$ \ are independent.
In the purely explosive case, the weak convergence in \eqref{help_rconv_2dnorm}
 was established in several other papers, see, e.g., Jeganathan \cite[Theorem 14]{Jeg}, Mikulski and Monsour \cite[Theorem 1]{MikMon1} or
 Monsour and Mikulski \cite[Theorem on page 146]{MikMon2}.

Since mixing convergence yields convergence in distribution,
 as a consequence of Theorem \ref{|lambda_1|>|lambda_2|}, we have
 \begin{equation}\label{|lambda_1|>|lambda_2|_distr}
   \begin{bmatrix}
     \left(\sum_{k=1}^n X_{k-1}^2\right)^{1/2}
      & \frac{\sum_{k=1}^n X_{k-1} X_{k-2}}{\left(\sum_{k=1}^n X_{k-1}^2\right)^{1/2}} \\[3mm]
     \frac{\sum_{k=1}^n X_{k-1} X_{k-2}}{\left(\sum_{k=1}^n X_{k-2}^2\right)^{1/2}}
      & \left(\sum_{k=1}^n X_{k-2}^2\right)^{1/2}
    \end{bmatrix}
    \begin{bmatrix}
     \hvartheta_1^{(n)} - \vartheta_1 \\
     \hvartheta_2^{(n)} - \vartheta_2
    \end{bmatrix} \\
   \distrP  N \begin{bmatrix} 1 \\ \sign(\lambda_1) \end{bmatrix} \qquad \text{as \ $n \to \infty$,}
 \end{equation}
 where \ $N$ \ is a standard normally distributed random variable.

Next we check that \eqref{help_rconv_2dnorm} yields \eqref{|lambda_1|>|lambda_2|_distr}.
  Using \eqref{LSE_diff_formula} and \eqref{langle_bM_rangle_n}, the weak convergence in \eqref{help_rconv_2dnorm} can be rewritten in the form
 \[
   \langle \bM\rangle_n^{1/2} \langle \bM\rangle_n^{-1} \bM_n
   = \langle \bM\rangle_n^{-1/2} \bM_n
   \distrP \cN_2(\bzero, \bI_2) \qquad \text{as \ $n \to \infty$.}
 \]
By the decomposition in Step 1 of the proof of Theorem \ref{|lambda_1|>|lambda_2|}, the weak convergence in \eqref{|lambda_1|>|lambda_2|_distr}
 can be rewritten in the form
 \begin{align}\label{help_AnMn}
  \bA_n \bM_n \distrP  N \begin{bmatrix} 1 \\ \sign(\lambda_1) \end{bmatrix} \qquad \text{as \ $n \to \infty$,}
 \end{align}
 where \ $\bA_n$, $n\in\NN$, \ is given in \eqref{help_An} with \ $\sigma=1$.
\ Since on the event \ $\Omega_n$ \ (given in \eqref{Omega_n}), we have
 \ $\bA_n \bM_n = \bA_n \langle \bM\rangle_n^{1/2} (\langle \bM\rangle_n^{-1/2}\bM_n)$, \
 to check that \eqref{help_rconv_2dnorm} yields \eqref{|lambda_1|>|lambda_2|_distr}, by Slutsky's lemma,
 it is enough to verify that
 \begin{align}\label{help_discussion}
   \bA_n \langle \bM\rangle_n^{1/2}
    \asP \frac{1}{\sqrt{1+\lambda_1^{-2}}} \begin{bmatrix}
                                             1 & \lambda_1^{-1} \\
                                             \sign(\lambda_1) & \vert\lambda_1\vert^{-1}  \\
                                           \end{bmatrix}
    \qquad \text{as \ $n\to\infty$.}
 \end{align}
Indeed,
 \begin{align*}
 & \left(\frac{1}{\sqrt{1+\lambda_1^{-2}}}
     \begin{bmatrix}
        1 & \lambda_1^{-1} \\
        \sign(\lambda_1) & \vert\lambda_1\vert^{-1}  \\
    \end{bmatrix} \right)
   \left(\frac{1}{\sqrt{1+\lambda_1^{-2}}}
     \begin{bmatrix}
        1 & \lambda_1^{-1} \\
        \sign(\lambda_1) & \vert\lambda_1\vert^{-1}  \\
    \end{bmatrix} \right)^\top \\
 & \qquad = \frac{1}{1+\lambda_1^{-2}}
     \begin{bmatrix}
        1 + \lambda_1^{-2} & \sign(\lambda_1)+(\lambda_1 \vert \lambda_1\vert)^{-1} \\
        \sign(\lambda_1)+(\lambda_1 \vert \lambda_1\vert)^{-1} & 1 + \lambda_1^{-2} \\
    \end{bmatrix}
  =  \begin{bmatrix}
        1 & \sign(\lambda_1) \\
        \sign(\lambda_1) & 1 \\
    \end{bmatrix} ,
 \end{align*}
 which coincides with the covariance matrix of \ $N [ 1 , \sign(\lambda_1) ]^\top$, \ as desired.
Next, we prove \eqref{help_discussion}.
Recall that if \ $\bV=(v_{i,j})_{i,j=1,2}\in\RR^{2\times 2}$ \ is a symmetric and positive definite matrix, then
 \[
  \bV^{1/2} = \frac{1}{\sqrt{ v_{1,1} + v_{2,2} + 2\sqrt{\det(\bV)}}}\Big(\bV + \sqrt{\det(\bV)} \bI_2\Big),
 \]
 and hence on the event \ $\Omega_n$, \
 \[
  \langle \bM\rangle_n^{1/2}
   = \frac{1}{\sqrt{ \sum_{k=1}^n X_{k-1}^2 + \sum_{k=1}^n X_{k-2}^2 + 2\sqrt{ \det(\langle \bM\rangle_n)}}}
      \Big( \langle \bM\rangle_n + \sqrt{ \det(\langle \bM\rangle_n) } \bI_2 \Big).
 \]
Consequently, by \eqref{sum_{k=1}^n X_{k-1}^2-}, \eqref{sum_{k=1}^n X_{k-2}^2-}, \eqref{bA_n}, \eqref{sum_{k=1}^n X_{k-1} X_{k-2}-},
 the absolute continuity of \ $Y$ \ (see Step 2 of Theorem \ref{|lambda_1|>|lambda_2|}) and Corollary \ref{Cor_detMn}, we have
 \begin{align*}
  &\bA_n \langle \bM\rangle_n^{1/2}
     = \frac{1}{\sqrt{ \sum_{k=1}^n X_{k-1}^2 +  \sum_{k=1}^n X_{k-2}^2 + 2\sqrt{ \det(\langle \bM\rangle_n)}}}\\
  &\phantom{\bA_n \langle \bM\rangle_n^{1/2}=\;}
          \times \left(  \begin{bmatrix}
               \left(\sum_{k=1}^n X_{k-1}^2\right)^{1/2}
                  & \frac{\sum_{k=1}^n X_{k-1} X_{k-2}}{\left(\sum_{k=1}^n X_{k-1}^2\right)^{1/2}} \\[3mm]
               \frac{\sum_{k=1}^n X_{k-1} X_{k-2}}{\left(\sum_{k=1}^n X_{k-2}^2\right)^{1/2}}
                  & \left(\sum_{k=1}^n X_{k-2}^2\right)^{1/2}
             \end{bmatrix}
             + \sqrt{ \det(\langle \bM\rangle_n) } \bA_n \right)\\
    & = \frac{1}{\sqrt{  \lambda_1^{-2n}\sum_{k=1}^n X_{k-1}^2 +  \lambda_1^{-2n} \sum_{k=1}^n X_{k-2}^2
                        + 2\sqrt{ \lambda_1^{-4n} \det(\langle \bM\rangle_n)}}}\\
    &\phantom{=\;}\times     \left(  \begin{bmatrix}
               \left(\lambda_1^{-2n}\sum_{k=1}^n X_{k-1}^2\right)^{1/2}
                  & \frac{\lambda_1^{-2n}\sum_{k=1}^n X_{k-1} X_{k-2}}{\left(\lambda_1^{-2n}\sum_{k=1}^n X_{k-1}^2\right)^{1/2}} \\[3mm]
               \frac{\lambda_1^{-2n}\sum_{k=1}^n X_{k-1} X_{k-2}}{\left(\lambda_1^{-2n}\sum_{k=1}^n X_{k-2}^2\right)^{1/2}}
                  & \left(\lambda_1^{-2n}\sum_{k=1}^n X_{k-2}^2\right)^{1/2}
             \end{bmatrix}
             + \sqrt{ \lambda_1^{-4n} \det(\langle \bM\rangle_n) } \lambda_1^n\bA_n \right)\\
     &\asP
       \frac{1}{\sqrt{\frac{1}{\lambda_1^2-1}Y^2 + \frac{1}{\lambda_1^2(\lambda_1^2-1)}Y^2 } }
           \begin{bmatrix}
             \frac{\vert Y\vert}{\sqrt{\lambda_1^2-1}} & \frac{Y^2/(\lambda_1(\lambda_1^2-1))}{\vert Y\vert/\sqrt{\lambda_1^2-1}}  \\[2mm]
             \frac{Y^2/(\lambda_1(\lambda_1^2-1))}{\vert Y\vert/(\vert \lambda_1\vert\sqrt{\lambda_1^2-1})}  &  \frac{\vert Y\vert}{\vert \lambda_1\vert\sqrt{\lambda_1^2-1}} \\
           \end{bmatrix}\\
     &\phantom{\asP} =  \frac{1}{\sqrt{1+\lambda_1^{-2}}} \begin{bmatrix}
                                             1 & \lambda_1^{-1} \\
                                             \sign(\lambda_1) & \vert\lambda_1\vert^{-1}  \\
                                           \end{bmatrix} \qquad \text{as \ $n\to\infty$,}
 \end{align*}
 yielding \eqref{help_discussion}, as desired.
Note that \eqref{help_discussion} holds in case of \ $\vert \lambda_2\vert=1$ \ as well
 (indeed, in the previous argument we did not use that \ $\vert \lambda_2\vert\ne 1$).

Now let us suppose that the conditions of Theorem \ref{|lambda_1|>|lambda_2|} hold together with \ $\vert \lambda_2\vert = 1$, \ i.e.,
 \ $\lambda_1\in\RR$, \ $\vert \lambda_1\vert>1$, \ and either \ $\lambda_2=1$ \ or \ $\lambda_2=-1$.
\ For simplicity, let us suppose again that \ $\sigma=1$.
Provided that the proof of part (c) of Theorem 3 in Monsour \cite{Mon} is complete, in our considered case it states that
   \begin{align}\label{help_rconv_Mon_limit}
   \begin{bmatrix}
     \sum_{k=1}^n X_{k-1}^2 & \sum_{k=1}^n X_{k-1} X_{k-2} \\[1mm]
     \sum_{k=1}^n X_{k-1} X_{k-2} & \sum_{k=1}^n X_{k-2}^2
    \end{bmatrix}^{1/2}
    \begin{bmatrix}
     \hvartheta_1^{(n)} - \vartheta_1 \\
     \hvartheta_2^{(n)} - \vartheta_2
    \end{bmatrix}
   \distrP  \cO^{(\pm)} \begin{bmatrix}
                 \zeta^{(\pm)} \\[1mm]
                 \frac{\pm\int_0^1 W^{(\pm)}_u\dd W^{(\pm)}_u}{(\int_0^1 (W^{(\pm)}_u)^2\,\dd u)^{1/2}} \\
               \end{bmatrix}
 \end{align}
 as \ $n \to \infty$, \ where the \ $\pm$ \ sign is according to \ $\lambda_2=\pm 1$, \
 $\cO^{(\pm)}$ \ is some random \ $2\times 2$ \ orthogonal matrix, \ $\zeta^{(\pm)}$ \ is a one-dimensional standard normally distributed random variable,
  \ $(W^{(\pm)}_u)_{u\in[0,1]}$ \ is standard Wiener process such that \ $\zeta^{(\pm)}$ \ and \ $\int_0^1 W^{(\pm)}_u\dd W^{(\pm)}_u$ \ are uncorrelated.
Here \ $\cO^{(\pm)}$ and \ $(\zeta^{(\pm)}, \int_0^1 W^{(\pm)}_u\dd W^{(\pm)}_u/(\int_0^1 (W^{(\pm)}_u)^2\,\dd u)^{1/2})$ \ are not necessarily independent.
Using \eqref{help_discussion} (which holds in case of \ $\vert \lambda_2\vert=1$ \ as well)
 and Slutsky's lemma, similarly as we have seen before, \eqref{help_rconv_Mon_limit} should yield that
 \begin{align*}
  &\begin{bmatrix}
     \left(\sum_{k=1}^n X_{k-1}^2\right)^{1/2}
      & \frac{\sum_{k=1}^n X_{k-1} X_{k-2}}{\left(\sum_{k=1}^n X_{k-1}^2\right)^{1/2}} \\[3mm]
     \frac{\sum_{k=1}^n X_{k-1} X_{k-2}}{\left(\sum_{k=1}^n X_{k-2}^2\right)^{1/2}}
      & \left(\sum_{k=1}^n X_{k-2}^2\right)^{1/2}
    \end{bmatrix}
    \begin{bmatrix}
     \hvartheta_1^{(n)} - \vartheta_1 \\
     \hvartheta_2^{(n)} - \vartheta_2
    \end{bmatrix}
    = \bA_n \bM_n
    = \bA_n \langle \bM\rangle_n^{\frac{1}{2}} \langle \bM\rangle_n^{-\frac{1}{2}}\bM_n
    \distrP\\
  &\qquad \distrP
     \frac{1}{\sqrt{1+\lambda_1^{-2}}} \begin{bmatrix}
                                             1 & \lambda_1^{-1} \\
                                             \sign(\lambda_1) & \vert\lambda_1\vert^{-1}  \\
                                        \end{bmatrix}
     \cO^{(\pm)}  \begin{bmatrix}
                 \zeta^{(\pm)} \\[1mm]
                 \frac{\pm\int_0^1 W^{(\pm)}_u\dd W^{(\pm)}_u}{(\int_0^1 (W^{(\pm)}_u)^2\,\dd u)^{1/2}} \\
          \end{bmatrix}     \qquad \text{as \ $n\to\infty$.}
 \end{align*}
Taking into account \eqref{conv_|lambda_1|>|lambda_2|} and that mixing convergence yields convergence in distribution,
 under the conditions of Theorem \ref{|lambda_1|>|lambda_2|} together with \ $\vert \lambda_2\vert = 1$ \ and \ $\sigma=1$,
 \ it should hold that
 \[
   \frac{1}{\sqrt{1+\lambda_1^{-2}}} \begin{bmatrix}
                                             1 & \lambda_1^{-1} \\
                                             \sign(\lambda_1) & \vert\lambda_1\vert^{-1}  \\
                                        \end{bmatrix}
     \cO^{(\pm)}  \begin{bmatrix}
                 \zeta^{(\pm)} \\[1mm]
                 \frac{\pm\int_0^1 W^{(\pm)}_u\dd W^{(\pm)}_u}{(\int_0^1 (W^{(\pm)}_u)^2\,\dd u)^{1/2}} \\
          \end{bmatrix}
         \distre  N \begin{bmatrix} 1 \\ \sign(\lambda_1) \end{bmatrix}.
 \]
We were not able to check whether the previous equality in distribution holds or not
 (mainly due to the lack of an explicit form for the random matrix \ $\cO^{(\pm)}$), \ and in fact, we are not sure that it is true, since,
 as we detailed earlier, we are not convinced that the proof of part (c) of Theorem 3 in Monsour \cite{Mon} is complete/correct.

To finish the comparison of our results and part (c) of Theorem 3 in Monsour \cite{Mon},
  we emphasize that both in \eqref{help_rconv_2dnorm} and in \eqref{|lambda_1|>|lambda_2|_distr} the type of convergence
 is convergence in distribution, and in \eqref{|lambda_1|>|lambda_2|} we proved mixing convergence which is stronger than convergence in distribution, so in general
 \eqref{help_rconv_2dnorm} (or \eqref{|lambda_1|>|lambda_2|_distr}) would not yield \eqref{|lambda_1|>|lambda_2|} without any additional work.
Furthermore, our proof technique is completely different from that of Monsour \cite{Mon}.

Now we turn to compare Theorem \ref{|lambda_1|>|lambda_2|} (of the present paper) and Th\'eor\`{e}me 1 in Touati \cite{Tou}
 specialized to the two-dimensional case.
Let us suppose that the conditions of Theorem \ref{|lambda_1|>|lambda_2|} hold together with
  \ $\vert \lambda_1\vert > \vert \lambda_2\vert>1$ \ (i.e., we consider the purely explosive case).
For simplicity, let us suppose that \ $\sigma=1$.
Th\'eor\`{e}me 1 in Touati \cite{Tou} specialized to this case states that
 \begin{align}\label{help_Touati_1}
  \bvartheta^{-n}
   \begin{bmatrix}
     \sum_{k=1}^n X_{k-1}^2 & \sum_{k=1}^n X_{k-1} X_{k-2} \\[1mm]
     \sum_{k=1}^n X_{k-1} X_{k-2} & \sum_{k=1}^n X_{k-2}^2
    \end{bmatrix}
    \begin{bmatrix}
     \hvartheta_1^{(n)} - \vartheta_1 \\
     \hvartheta_2^{(n)} - \vartheta_2
    \end{bmatrix}
   \distrP  \sum_{k=1}^\infty \bvartheta^{-k} \widetilde\bY\zeta_k \qquad \text{as \ $n \to \infty$,}
 \end{align}
 where
 \[
 \widetilde\bY
    := \begin{bmatrix}
         X_0 \\
         X_{-1} \\
       \end{bmatrix}
       + \sum_{k=1}^\infty \bvartheta^{-k}
             \begin{bmatrix}
               Z_k \\
               0 \\
             \end{bmatrix},
 \]
 and \ $(\zeta_k)_{k\in\NN}$ \ is a sequence of independent and standard normally distributed random variables
 such that \ $(\zeta_k)_{k\in\NN}$ \ is independent of \ $(X_0,X_{-1})$ \ and \ $(Z_n)_{n\in\NN}$ \
 (and consequently, \ $(\zeta_k)_{k\in\NN}$ \ is independent of \ $\widetilde\bY$ \ as well).
Recall that, since mixing convergence yields convergence in distribution,
 as a consequence of Theorem \ref{|lambda_1|>|lambda_2|}, we have \eqref{|lambda_1|>|lambda_2|_distr}.

Next, we check that \eqref{help_Touati_1} yields \eqref{|lambda_1|>|lambda_2|_distr} in the considered purely explosive case.
Using \eqref{LSE_diff_formula} and \eqref{langle_bM_rangle_n}, the weak convergence in \eqref{help_Touati_1} can be rewritten in the form
 \begin{align}\label{help_Touati_2}
   \bvartheta^{-n} \langle \bM\rangle_n (\langle \bM\rangle_n^{-1} \bM_n)
      = \bvartheta^{-n} \bM_n
      \distrP \sum_{k=1}^\infty \bvartheta^{-k} \widetilde\bY\zeta_k  \qquad \text{as \ $n \to \infty$.}
 \end{align}
We have already seen that the weak convergence in \eqref{|lambda_1|>|lambda_2|_distr}
 can be rewritten in the form \eqref{help_AnMn}.
Let us consider the decomposition \ $\bA_n \bM_n = \bA_n \bvartheta^{n} (\bvartheta^{-n}\bM_n)$, \ $n\in\NN$.
\ By \eqref{bvartheta^n}, \eqref{bA_n} and using that \ $\big(\frac{\lambda_2}{\lambda_1}\big)^n\to 0$ \ as \ $n\to\infty$ \
 (due to \ $\vert \lambda_2\vert < \vert \lambda_1\vert$), \ we have
 \begin{align*}
 \bA_n \bvartheta^{n}
    & = \lambda_1^n\bA_n\frac{1}{\lambda_1-\lambda_2}
      \begin{bmatrix}
       \lambda_1 & - \lambda_1 \lambda_2 \\
       1 & - \lambda_2
      \end{bmatrix}
      + \left(\frac{\lambda_2}{\lambda_1}\right)^n\lambda_1^n \bA_n\frac{1}{\lambda_1-\lambda_2}
        \begin{bmatrix}
         - \lambda_2 & \lambda_1 \lambda_2 \\
         - 1 & \lambda_1
        \end{bmatrix} \\
    & \asP
      \frac{\sqrt{\lambda_1^2-1}}{|Y|}
      \begin{bmatrix}
       1 & 0 \\
       0 & |\lambda_1|
      \end{bmatrix}
      \frac{1}{\lambda_1-\lambda_2}
      \begin{bmatrix}
       \lambda_1 & - \lambda_1 \lambda_2 \\
       1 & - \lambda_2
      \end{bmatrix}
      \qquad \text{as \ $n\to\infty$,}
 \end{align*}
 where \ $Y$ \ is given in \eqref{conv_X_n}.
Hence, by Slutsky's lemma, \eqref{help_Touati_2} yields that
 \begin{align*}
   \bA_n \bM_n
     \distrP
     \frac{\sqrt{\lambda_1^2-1}}{|Y|(\lambda_1-\lambda_2)}
      \begin{bmatrix}
       1 & 0 \\
       0 & |\lambda_1|
      \end{bmatrix}
      \begin{bmatrix}
       \lambda_1 & - \lambda_1 \lambda_2 \\
       1 & - \lambda_2
      \end{bmatrix}
      \sum_{k=1}^\infty \bvartheta^{-k} \widetilde\bY\zeta_k
      \qquad \text{as \ $n\to\infty$.}
 \end{align*}
Using \eqref{bvartheta^n}, for each \ $k\in\NN$, \ we have
 \begin{align*}
  \bvartheta^{-k}
   &= \begin{bmatrix}
       \lambda_+ & \lambda_- \\
       1 & 1
      \end{bmatrix}
      \begin{bmatrix}
       \lambda_+^{-k} & 0 \\
       0 & \lambda_-^{-k}
      \end{bmatrix}
      \begin{bmatrix}
       \lambda_+ & \lambda_- \\
       1 & 1
      \end{bmatrix}^{-1} \\
   &= \frac{\lambda_1^{-k}}{\lambda_1-\lambda_2}
      \begin{bmatrix}
       \lambda_1 & - \lambda_1 \lambda_2 \\
       1 & - \lambda_2
      \end{bmatrix}
      + \frac{\lambda_2^{-k}}{\lambda_1-\lambda_2}
        \begin{bmatrix}
         - \lambda_2 & \lambda_1 \lambda_2 \\
         - 1 & \lambda_1
        \end{bmatrix}.
 \end{align*}
Hence for each \ $k\in\NN$, \ we get
 \begin{align}\label{help_vartheta_inverz}
  \begin{split}
     \begin{bmatrix}
       \lambda_1 & - \lambda_1 \lambda_2 \\
       1 & - \lambda_2
      \end{bmatrix}
      \bvartheta^{-k}
    &= \frac{\lambda_1^{-k}}{\lambda_1-\lambda_2}
      \begin{bmatrix}
       \lambda_1 & - \lambda_1 \lambda_2 \\
       1 & - \lambda_2
      \end{bmatrix}^2
      + \frac{\lambda_2^{-k}}{\lambda_1-\lambda_2}
        \begin{bmatrix}
       \lambda_1 & - \lambda_1 \lambda_2 \\
       1 & - \lambda_2
      \end{bmatrix}
        \begin{bmatrix}
         - \lambda_2 & \lambda_1 \lambda_2 \\
         - 1 & \lambda_1
        \end{bmatrix}   \\
      &= \frac{\lambda_1^{-k}}{\lambda_1-\lambda_2}
        \begin{bmatrix}
       \lambda_1 & - \lambda_1 \lambda_2 \\
       1 & - \lambda_2
      \end{bmatrix}^2.
  \end{split}
 \end{align}
Consequently,
 \begin{align*}
   \bA_n \bM_n
     \distrP
     \frac{\sqrt{\lambda_1^2-1}}{|Y|(\lambda_1-\lambda_2)}
      \begin{bmatrix}
       1 & 0 \\
       0 & |\lambda_1|
      \end{bmatrix}
      \sum_{k=1}^\infty
        \frac{\lambda_1^{-k}}{\lambda_1-\lambda_2}
        \begin{bmatrix}
           \lambda_1 & - \lambda_1 \lambda_2 \\
             1 & - \lambda_2
         \end{bmatrix}^2 \widetilde\bY\zeta_k
        \qquad \text{as \ $n\to\infty$.}
 \end{align*}
Here, using \eqref{help_vartheta_inverz}, we have
 \begin{align*}
   \begin{bmatrix}
           \lambda_1 & - \lambda_1 \lambda_2 \\
             1 & - \lambda_2
   \end{bmatrix} \widetilde\bY
  & =
   \begin{bmatrix}
           \lambda_1 & - \lambda_1 \lambda_2 \\
             1 & - \lambda_2
   \end{bmatrix}
    \left(\begin{bmatrix}
         X_0 \\
         X_{-1} \\
       \end{bmatrix}
       + \sum_{\ell=1}^\infty \bvartheta^{-\ell}
             \begin{bmatrix}
               Z_\ell \\
               0 \\
             \end{bmatrix}
   \right)\\
  & = \begin{bmatrix}
           \lambda_1 & - \lambda_1 \lambda_2 \\
             1 & - \lambda_2
   \end{bmatrix}
   \begin{bmatrix}
         X_0 \\
         X_{-1} \\
   \end{bmatrix}
   + \sum_{\ell=1}^\infty \frac{\lambda_1^{-\ell}}{\lambda_1-\lambda_2}
             \begin{bmatrix}
                \lambda_1 & - \lambda_1 \lambda_2 \\
                1 & - \lambda_2
             \end{bmatrix}^2
             \begin{bmatrix}
               Z_\ell \\
               0 \\
             \end{bmatrix}\\
  & = \begin{bmatrix}
           \lambda_1 & - \lambda_1 \lambda_2 \\
             1 & - \lambda_2
       \end{bmatrix}
       \left(
          \begin{bmatrix}
            X_0 \\
            X_{-1} \\
          \end{bmatrix}
          + \frac{1}{\lambda_1-\lambda_2} \left(\sum_{\ell=1}^\infty \lambda_1^{-\ell} Z_\ell \right)
            \begin{bmatrix}
            \lambda_1 \\
            1 \\
          \end{bmatrix}
       \right).
 \end{align*}
This yields that
 \begin{align*}
 \bA_n \bM_n
     \distrP
      &\,\frac{\sqrt{\lambda_1^2-1}}{|Y|(\lambda_1-\lambda_2)}
      \begin{bmatrix}
       1 & 0 \\
       0 & |\lambda_1|
      \end{bmatrix}
       \begin{bmatrix}
           \lambda_1 & - \lambda_1 \lambda_2 \\
             1 & - \lambda_2
       \end{bmatrix}^{2}\\
     &\times \left(
        \begin{bmatrix}
            X_0 \\
            X_{-1} \\
          \end{bmatrix}
          + \frac{1}{\lambda_1-\lambda_2} \left(\sum_{\ell=1}^\infty \lambda_1^{-\ell} Z_\ell \right)
            \begin{bmatrix}
            \lambda_1 \\
            1 \\
          \end{bmatrix}
       \right)
       \sum_{k=1}^\infty \frac{\lambda_1^{-k}}{\lambda_1-\lambda_2}\zeta_k
 \end{align*}
 as \ $n\to\infty$.
\ Here
 \begin{align*}
  \begin{bmatrix}
       1 & 0 \\
       0 & |\lambda_1|
  \end{bmatrix}
  &\begin{bmatrix}
   \lambda_1 & - \lambda_1 \lambda_2 \\
   1 & - \lambda_2
  \end{bmatrix}^{2}
  \left(
        \begin{bmatrix}
            X_0 \\
            X_{-1} \\
          \end{bmatrix}
          + \frac{1}{\lambda_1-\lambda_2} \left(\sum_{\ell=1}^\infty \lambda_1^{-\ell} Z_\ell \right)
            \begin{bmatrix}
            \lambda_1 \\
            1 \\
          \end{bmatrix}
  \right) \\
 &= \begin{bmatrix}
       1 & 0 \\
       0 & |\lambda_1|
  \end{bmatrix}
  \begin{bmatrix}
   \lambda_1 & - \lambda_1 \lambda_2 \\
   1 & - \lambda_2
  \end{bmatrix}
  \left(
   \begin{bmatrix}
       \lambda_1(X_0 - \lambda_2 X_{-1}) \\
       X_0 - \lambda_2 X_{-1}
  \end{bmatrix}
  + \left( \sum_{\ell=1}^\infty \lambda_1^{-\ell} Z_\ell \right)
   \begin{bmatrix}
       \lambda_1 \\
        1
    \end{bmatrix}
  \right)\\
 &= \begin{bmatrix}
       1 & 0 \\
       0 & |\lambda_1|
  \end{bmatrix}
  \begin{bmatrix}
   \lambda_1 & - \lambda_1 \lambda_2 \\
   1 & - \lambda_2
  \end{bmatrix}
  \begin{bmatrix}
    (\lambda_1 - \lambda_2) Y \\
    \frac{\lambda_1 - \lambda_2}{\lambda_1} Y \\
  \end{bmatrix}\\
  &=
   (\lambda_1 - \lambda_2)\begin{bmatrix}
       1 & 0 \\
       0 & |\lambda_1|
  \end{bmatrix}
  \begin{bmatrix}
   \lambda_1 & - \lambda_1 \lambda_2 \\
   1 & - \lambda_2
  \end{bmatrix}
  \begin{bmatrix}
    1 \\
    \frac{1}{\lambda_1} \\
  \end{bmatrix}
  Y.
 \end{align*}
Consequently,
 \begin{align*}
  \bA_n \bM_n
    \distrP  \frac{\sqrt{\lambda_1^2-1}}{\lambda_1-\lambda_2}
             \frac{Y}{\vert Y\vert}
             \begin{bmatrix}
                1 & 0 \\
                0 & |\lambda_1|
             \end{bmatrix}
             \begin{bmatrix}
                \lambda_1 & - \lambda_1 \lambda_2 \\
                1 & - \lambda_2
             \end{bmatrix}
             \begin{bmatrix}
                 1 \\
              \frac{1}{\lambda_1} \\
             \end{bmatrix}
             \sum_{k=1}^\infty \lambda_1^{-k}\zeta_k
    \qquad \text{as \ $n\to\infty$.}
 \end{align*}
Recall that \ $(\zeta_k)_{k\in\NN}$ \ is independent of \ $(X_0,X_{-1})$ \ and \ $(Z_n)_{n\in\NN}$, \ and hence
 \ $(\zeta_k)_{k\in\NN}$ \ is independent of \ $Y$ \ as well.
Consequently, since \ $ \frac{Y}{\vert Y\vert} = \sign(Y)$, \ $\sum_{k=1}^\infty \lambda_1^{-k}\zeta_k$ \ is normally distributed with mean zero,
 and the law of a normally distributed random variable with mean zero is invariant under an orthogonal transformation, we have
 \begin{align*}
  \bA_n \bM_n
    \distrP
     &\; \frac{\sqrt{\lambda_1^2-1}}{\lambda_1-\lambda_2}
      \begin{bmatrix}
                1 & 0 \\
                0 & |\lambda_1|
      \end{bmatrix}
      \begin{bmatrix}
                \lambda_1 - \lambda_2 \\
                1 - \frac{\lambda_2}{\lambda_1}
      \end{bmatrix}
      \sum_{k=1}^\infty \lambda_1^{-k}\zeta_k\\
     &=\sqrt{\lambda_1^2-1}
       \sum_{k=1}^\infty \lambda_1^{-k}\zeta_k
       \begin{bmatrix}
                 1 \\
              \sign(\lambda_1) \\
       \end{bmatrix}
      \qquad \text{as \ $n\to\infty$,}
 \end{align*}
 where \ $\sqrt{\lambda_1^2-1}\sum_{k=1}^\infty \lambda_1^{-k}\zeta_k$ \ is normally distributed with mean zero and variance
 \begin{align*}
   (\lambda_1^2-1)\sum_{k=1}^\infty \lambda_1^{-2k}
     = (\lambda_1^2-1) \frac{\lambda_1^{-2}}{1-\lambda_1^{-2}}
     =1,
 \end{align*}
 i.e., \ $\sqrt{\lambda_1^2-1}\sum_{k=1}^\infty \lambda_1^{-k}\zeta_k$ \ is standard normally distributed.
This yields \eqref{help_AnMn}, as desired.

Next we compare our proof technique with that of Th\'eor\`{e}me 1 in Touati \cite{Tou}.
The proofs for the purely explosive case in Th\'eor\`{e}me 1 in Touati \cite{Tou} are based on a limit theorem for a triangular array
 which does not satisfy the asymptotically negligible condition stating stable convergence of the row-sums of the triangular array in question,
 see Touati \cite[Th\'eor\`{e}me A]{Tou}.
Our Theorem \ref{MSLTES}  might be considered to be similar to Th\'eor\`{e}me A in Touati \cite{Tou}, but in fact
 these two theorems are quite different.
For example, in Th\'eor\`{e}me A in Touati \cite{Tou} there is a condition on the conditional characteristic function
 of the row-sums of the triangular array, while in our  Theorem \ref{MSLTES} (specified to triangular arrays)
 one can find a similar condition only on the conditional charasteristic function of the summands.
In the proof of Th\'eor\`{e}me 1 in Touati \cite{Tou} the author effectively used his Th\'eor\`{e}me A for proving  convergence in distribution
 of the appropriately normalized LSE of the AR parameters, but stable (mixing) convergence of the LSE in question
 was not proved.
We also note that our random normalization is different from the ones \ $\bvartheta^{-n} \langle\bM\rangle_n$ \ and \ $\bvartheta^n$ \ in Touati \cite[Th\'eor\`{e}me 1]{Tou},
  which contain the unknown parameters to be estimated, while
 our random normalization in Theorem \ref{|lambda_1|>|lambda_2|} does not contain the unknown parameters, and hence
 one can easily use Theorem \ref{|lambda_1|>|lambda_2|} for constructing confidence regions, see Remark \ref{Rem_confidence}.
In the partially explosive case, according to Touati \cite[Section 3]{Tou}, an AR(2) process can be decomposed into a supercritical and a subcritical
 AR(1) process, and the results on the asymptotic behaviour of the LSE of the AR parameter for a supercritical and a subcritical
 AR(1) process were combined in Touati \cite[Section 3]{Tou}.
Hence our proof technique is completely different from that of Touati \cite[Section 3]{Tou} for partially explosive AR(2) processes.
Further, note that our Theorem \ref{|lambda_1|>|lambda_2|} also covers supercritical AR(2) processes with a unit root, while
 Touati \cite{Tou} did not handle this case.

Finally, we note that we tried to prove mixing (or stable-) convergence of \ $(\hvartheta_1^{(n)},\hvartheta_2^{(n)})$ \ as \ $n\to\infty$ \
 in the supercritical case using a non-random normalization, but our attempts have not been successful so far.
It is mainly due to the fact that the almost sure limit of \ $\lambda_1^{-2n}\langle \bM\rangle_n$ \ as \ $n\to\infty$ \ in Corollary \ref{Cor_detMn}
 is a random \ $\RR^{2\times 2}$-valued matrix having determinant zero \ $\PP$-almost surely.
We also note that, for our proof technique of Theorem \ref{|lambda_1|>|lambda_2|}, one needs the assumption that
 the characteristic roots of the AR(2) process in question are real and have different absolute values.
Indeed, for example, in formulae \eqref{bvartheta^n} and \eqref{help_new1} we divide by the difference \ $\lambda_1-\lambda_2$ \
 of the two characteristic roots, or in Step 2 (of the proof of Theorem \ref{|lambda_1|>|lambda_2|}) we also use that $\vert \lambda_2\vert/\vert \lambda_1\vert<1$.

\section{Simulation results}\label{section_sim}

In this section we illustrate Theorem \ref{|lambda_1|>|lambda_2|} using generated sample paths of a supercritical AR(2) process.
More precisely, by simulations, we illustrate the weak convergence
  \begin{equation}\label{|lambda_1|>|lambda_2|_distr_duplicate}
   \begin{bmatrix}
     \left(\sum_{k=1}^n X_{k-1}^2\right)^{1/2}
      & \frac{\sum_{k=1}^n X_{k-1} X_{k-2}}{\left(\sum_{k=1}^n X_{k-1}^2\right)^{1/2}} \\[3mm]
     \frac{\sum_{k=1}^n X_{k-1} X_{k-2}}{\left(\sum_{k=1}^n X_{k-2}^2\right)^{1/2}}
      & \left(\sum_{k=1}^n X_{k-2}^2\right)^{1/2}
    \end{bmatrix}
    \begin{bmatrix}
     \hvartheta_1^{(n)} - \vartheta_1 \\
     \hvartheta_2^{(n)} - \vartheta_2
    \end{bmatrix} \\
   \distrP  N \begin{bmatrix} 1 \\ \sign(\lambda_1) \end{bmatrix} \qquad \text{as \ $n \to \infty$,}
 \end{equation}
 where \ $N$ \ is a standard normally distributed random variable.
The weak convergence \eqref{|lambda_1|>|lambda_2|_distr_duplicate} is a consequence of the mixing convergence
 \eqref{conv_|lambda_1|>|lambda_2|} in Theorem \ref{|lambda_1|>|lambda_2|} with \ $\sigma=1$, \ since mixing convergence yields convergence in distribution.
We consider a supercritical AR(2) process \ $(X_n)_{n\geq -1}$ \  given in \eqref{AR(2)} such that \ $(Z_n)_{n\in\NN}$ \ is a sequence of
 i.i.d.\ random variables with \ $Z_1\distre \cN(0,1)$ \ and we suppose \ $(X_0,X_{-1})=(0,0)$.
\ We choose the AR parameters \ $\vartheta_1$ \ and \ $\vartheta_2$ \ in a way that the characteristic polynomial of the corresponding AR(2) process
 has real roots with different absolute values, and the AR(2) process in question is purely explosive (Case 1), partially explosive (Case 2), supercritical with
 a characteristic root \ $1$ \ (Case 3), and supercritical with a characteristic root \ $-1$ \ (Case 4), respectively.
In Table \ref{Table1} we summarize these four cases, and we give our particular choices of \ $(\vartheta_1,\vartheta_2)$ \ in these cases.

 \begin{table}[ht!]
\centering
\begin{small}
\begin{tabular}{|c|c|c|c|}
  \hline
      &  &  &  $(\vartheta_1,\vartheta_2)$   \\
  \hline
   Case 1: & $\vert \lambda_1\vert >1$ & $\vert\lambda_2\vert>1$  &   (1,3)  \\
  \hline
  Case 2: & $\vert\lambda_1\vert>1$ & $\vert\lambda_2\vert<1$ & (1,1)   \\
  \hline
  Case 3: & $\vert\lambda_1\vert>1$ & $\lambda_+=1$ \ or \ $\lambda_-=1$ & (-1,2)   \\
  \hline
  Case 4: & $\vert\lambda_1\vert>1$ & $\lambda_+=-1$ \ or \ $\lambda_-=-1$  & (-3,-2) \\
  \hline
\end{tabular}
\end{small}
\caption{The four cases with two real characteristic roots having different absolute values:
  purely explosive (Case 1), partially explosive (Case 1), supercritical with
         a characteristic root \ $1$ \ (Case 3), and supercritical with a characteristic root \ $-1$
         \ (Case 4). The fourth column contains our particular choices of \ $(\vartheta_1,\vartheta_2)$ \ in Cases 1--4. }
\label{Table1}
\end{table}

For every particular choice of \ $(\vartheta_1,\vartheta_2)$ \ we generate \ $1000$ \ independent $102$-length trajectories of
 \ $(X_n)_{n\geq -1}$, \ and for each of them we calculate the corresponding LSE \ $(\hvartheta_1^{(100)}, \hvartheta_2^{(100)})$ \
 and the left hand side of \eqref{|lambda_1|>|lambda_2|_distr_duplicate} with \ $n=100$, \ which we simply call
 the scaled error of the LSE in question.
We used the open software R for making the simulations.
Since the AR(2) process in question is supercritical, the generated value for \ $X_{100}$ \ is quite large (for our particular choices of
 parameters, its order is between \ $10^{20}$ \ and \ $10^{35}$),
 \ and, because of this, we used a precision of 800 bits for the calculations.
By making simulations, we realized that the two coordinates of the scaled error in question are almost the same
 (in case of \ $\lambda_1>1$) \ or are almost the $(-1)$-time of each other (in case of \ $\lambda_1<-1$)
 \ despite the high precision.
More precisely, for the 1000 generated trajectories, the Euclidean distance between the two vectors consisting of the first and second
 coordinate of the scaled error in question, respectively, is of order between \ $10^{-51}$ \ and \ $10^{-28}$.
\ As a consequence, the descriptive statistics and p-values of some tests for the scaled error of the LSE \ $\hvartheta_1^{(100)}$ \
 based on our 1000 replications and those for the scaled error of the LSE \ $\hvartheta_2^{(100)}$ \
 based on our 1000 replications are basically the same when we use the built-in functions of the software R.
Because of this, in what follows, we do not distinguish between the scaled errors of the LSE \ $\hvartheta_1^{(100)}$ \
 and that of the LSE \ $\hvartheta_2^{(100)}$, \ and we will simply write ''scaled errors of the LSEs of the parameters''.
There will be one exception when we calculate the empirical covariance between the scaled errors of the LSE
 \ $\hvartheta_1^{(100)}$ \ and that of the LSE \ $\hvartheta_2^{(100)}$.

For every particular choice of \ $(\vartheta_1,\vartheta_2)$, \
 we plot the density histogram of the scaled errors of the LSEs of the parameters based on the 1000 replications,
 and we also plot the density function of a standard normally distributed random variable in red on the same figure.
We calculate the empirical mean of the LSEs \ $(\hvartheta_1^{(100)},\hvartheta_2^{(100)})$ \ based on our 1000 replications.
We also calculate the empirical mean, variance, median, skewness, kurtosis and interquartile range of the scaled errors of the LSEs
 of the parameters based on the 1000 replications.
These values are expected to be around \ $0$, \ $1$, \ $0$, \ $0$, \ $3$ \ and \ $\Phi^{-1}(0.75) - \Phi^{-1}(0.25)\approx 1.349$ \
 (where \ $\Phi$ \ denotes the distribution function of a standard normally distributed random variable),
 since, due to \eqref{|lambda_1|>|lambda_2|_distr_duplicate}, the scaled errors of the LSEs of the parameters converge
 in distribution to a standard normal distribution as \ $n\to\infty$, \ and the (theoretical) mean, variance, median,
 skewness, kurtosis, and interquartile range of a standard normal distribution are the above listed values.
We also calculate the empirical covariance between the scaled errors of the LSEs \ $\hvartheta_1^{(100)}$ \ and \ $\hvartheta_2^{(100)}$,
 \ which is expected to be around \ $1$ \ if \ $\lambda_1>1$ \ and around \ $-1$ \ if \ $\lambda_1<-1$ \ due to \eqref{|lambda_1|>|lambda_2|_distr_duplicate}.
Using Kolmogorov-Smirnov, Pearson's chi-squared, Anderson-Darling and Jarque-Bera tests, we test whether the scaled errors of the LSEs
 of the parameters follow a normal distribution based on the 1000 replications of the \ $102$-length trajectories
 of the AR(2) process in question.
Using Kolmogorov-Smirnov test, we also test whether the scaled errors of the LSEs of the parameters
 follow a standard normal distribution.
We collect the \ $p$-values for these tests in a corresponding table.

For each of our four particular choices of \ $(\vartheta_1,\vartheta_2)$ \ (given in Table \ref{Table1}), we calculate the corresponding values of
 \ $\lambda_1$ \ and \ $\lambda_2$, \ and the empirical mean of the LSEs \ $(\hvartheta_1^{(100)},\hvartheta_2^{(100)})$ \ based on our 1000 replications,
 see Table \ref{becsles}.
 \begin{table}[ht!]
\centering
\begin{small}
\begin{tabular}{|c|c|c|c|}
  \hline
         & $\lambda_1$ & $\lambda_2$ &  Empirical mean of the LSEs \ $(\hvartheta_1^{(100)},\hvartheta_2^{(100)})$ \\
  \hline
   $(\vartheta_1,\vartheta_2) = (1,3)$  & 2.302776  & -1.302776 &  (0.999999,\;3.000001) \\
  \hline
   $(\vartheta_1,\vartheta_2) = (1,1)$  & 1.618034 & -0.618034 & (1.013616,\;0.977967) \\
  \hline
  $(\vartheta_1,\vartheta_2) = (-1,2)$  & -2 & 1 & (-1.016926,\;1.966146)\\
  \hline
  $(\vartheta_1,\vartheta_2) = (-3,-2)$  & -2 & -1 & (-2.980015,\;-1.960030) \\
  \hline
\end{tabular}
\end{small}
\caption{Values of $\lambda_1$ \ and \ $\lambda_2$, \ and the empirical mean of the LSEs \ $(\hvartheta_1^{(100)},\hvartheta_2^{(100)})$ \ based on 1000 replications.}
\label{becsles}
\end{table}

\begin{figure}[ht!]
\centering

\begin{subfigure}[t]{0.3\textwidth}
\centering
\includegraphics[width=\textwidth,height=5cm]{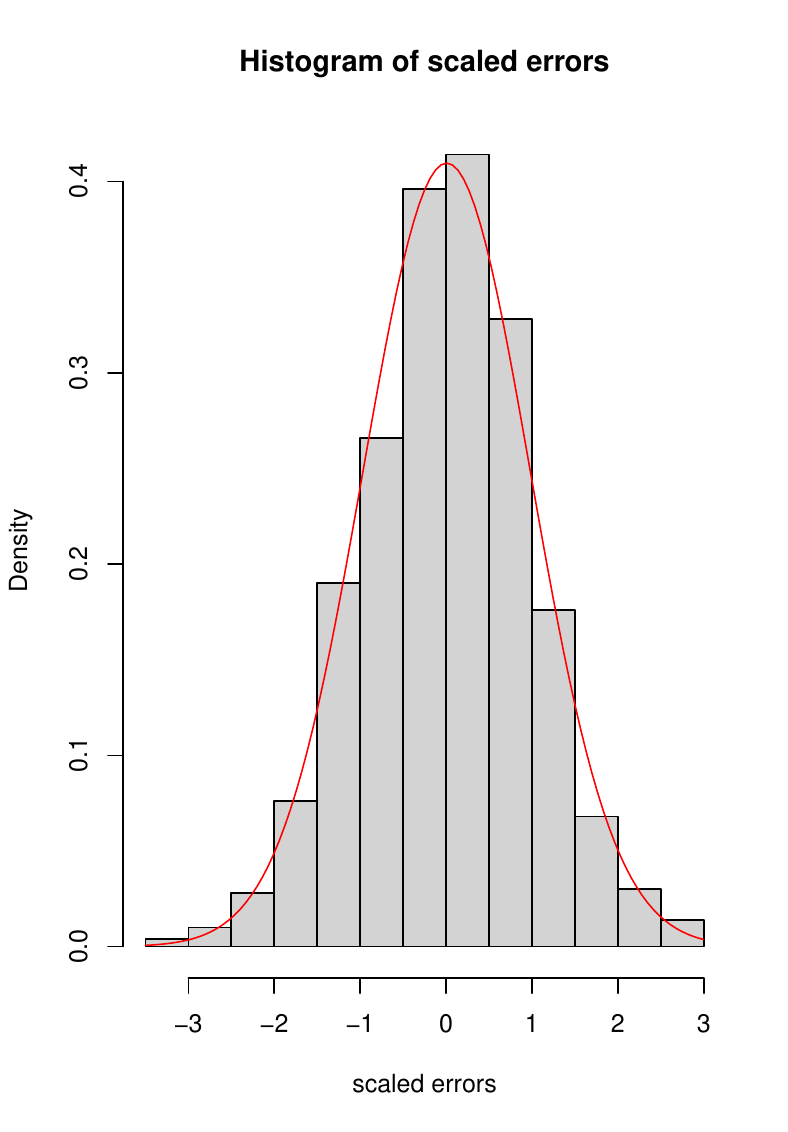}
\caption{$(\vartheta_1,\vartheta_2) = (1,3)$}
\end{subfigure}%
\hspace{2cm}
\begin{subfigure}[t]{0.3\textwidth}
\centering
\includegraphics[width=\textwidth,height=5cm]{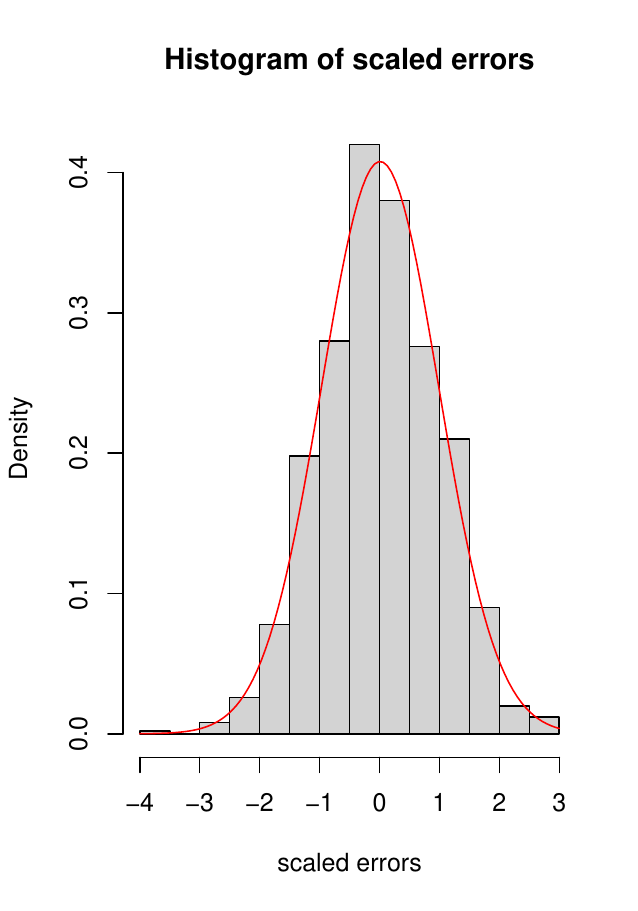}
\caption{$(\vartheta_1,\vartheta_2) = (1,1)$}
\end{subfigure}

\bigskip

\begin{subfigure}[t]{0.3\textwidth}
\centering
\includegraphics[width=\textwidth,height=5cm]{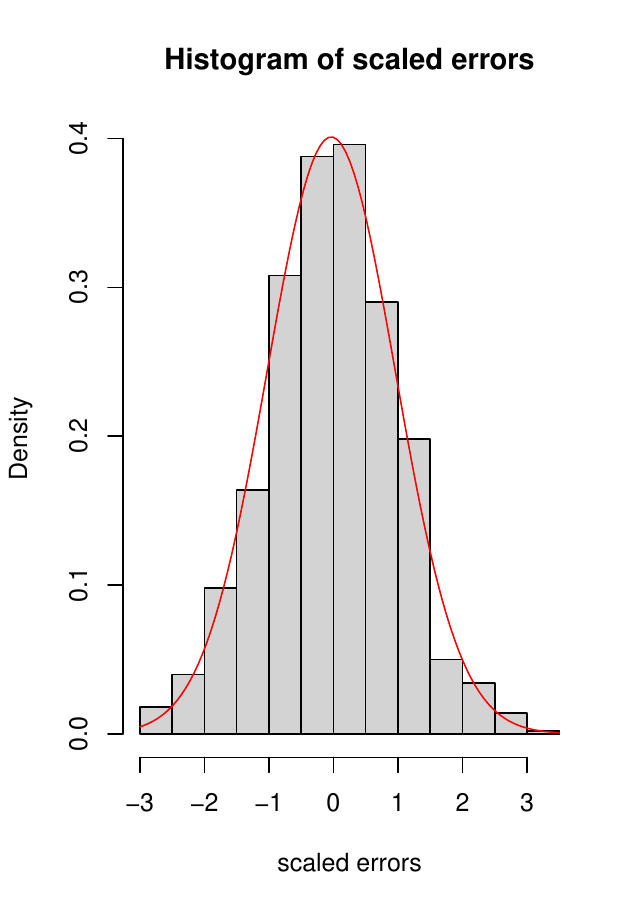}
\caption{$(\vartheta_1,\vartheta_2) = (-1,2)$}
\end{subfigure}%
\hspace{2cm}
\begin{subfigure}[t]{0.3\textwidth}
\centering
\includegraphics[width=\textwidth,height=5cm]{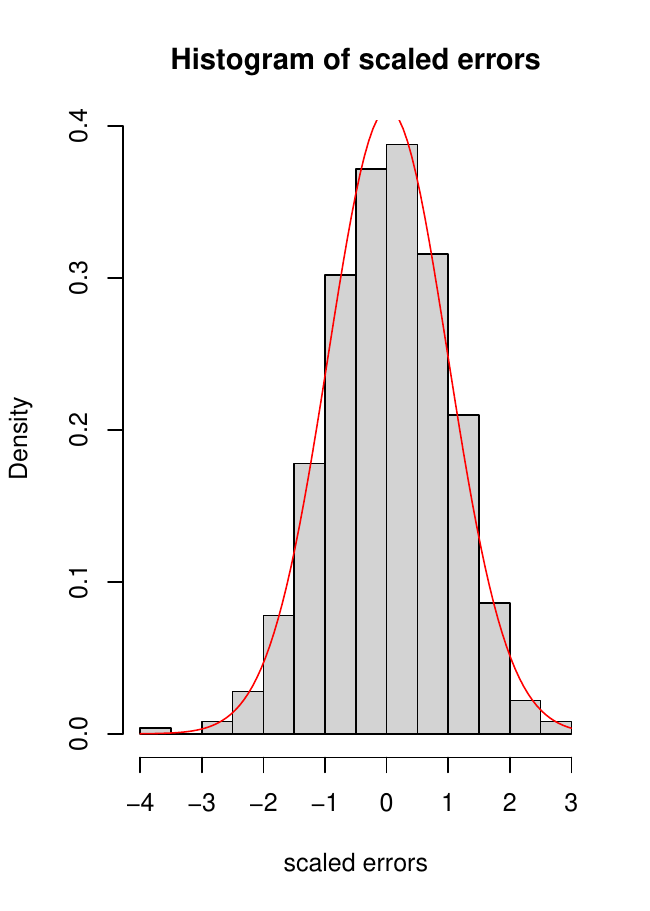}
\caption{$(\vartheta_1,\vartheta_2) = (-3,-2)$}
\end{subfigure}

\caption{Density histogram for the scaled errors.}
 \label{fig_hist}
\end{figure}

Figure \ref{fig_hist} contains the density histograms for the scaled errors for our four particular choices of
 \ $(\vartheta_1,\vartheta_2)$, \ and it supports the standard normality
 of the scaled errors of the LSEs of the parameters for each of the four particular choices of \ $(\vartheta_1,\vartheta_2)$.

Table \ref{descriptive_stats} contains the empirical mean, variance, median, skewness, kurtosis, interquartile range (IQR) and covariance for the scaled errors
 of the parameters for each of the four particular choices of \ $(\vartheta_1,\vartheta_2)$, \ and, in order to help the readers, in the second row of the table
 we present the corresponding theoretical values as well.
The simulation results in Table \ref{descriptive_stats} support the standard normality of the scaled errors
 of the LSEs of the considered four parameters.

 \begin{table}[ht!]
\centering
\begin{small}
\begin{tabular}{|c|c|c|c|c|c|c|c|}
  \hline
  Descriptive statistics  & Mean & Variance & Median & Skewness & Kurtosis & IQR & Covariance\\
  \hline
  Theoretical values & 0 & 1 & 0 & 0 & 3 & $\approx 1.349$  & $\sign(\lambda_1)$ \\
  \hline
   $(\vartheta_1,\vartheta_2) = (1,3)$  & 0.00753  &  0.94871 & 0.04073 &  -0.06399 & 3.03043 & 1.30492 & 0.94776 \\
  \hline
   $(\vartheta_1,\vartheta_2) = (1,1)$  & 0.01016 &  0.95664 & -0.01613  & -0.03674 & 2.97045 & 1.30881  &  0.95568  \\
  \hline
  $(\vartheta_1,\vartheta_2) = (-1,2)$  & -0.03357  &  0.98938 & -0.02913  & -0.03199  & 3.16162  & 1.30595 & -0.98839 \\
  \hline
  $(\vartheta_1,\vartheta_2) = (-3,-2)$  & 0.02392 & 0.94096  &  0.02936  & -0.12024 & 3.00297 & 1.31009  &  -0.94003 \\
  \hline
\end{tabular}
\end{small}
\caption{Empirical mean, variance, median, skewness, kurtosis, interquartile range and covariance for the scaled errors.
 The second row contains the corresponding theoretical values. }
\label{descriptive_stats}
\end{table}

Table \ref{tests} contains the \ $p$-values of the tests mentioned earlier for our four particular choices of
 \ $(\vartheta_1,\vartheta_2)$ \ (given in Table \ref{Table1}),
 and it shows that at any reasonable significance level all the considered four tests accept normality,
 and the Kolmogorov-Smirnov test accepts standard normality of the scaled errors in question.

 \begin{table}[ht!]
\centering
\begin{small}
\begin{tabular}{|c|c|c|c|c|c|}
  \hline
  $p$-values of normality tests  & KS test & KS test for $\cN(0,1)$ & PCS test & AD test & JB test\\
  \hline
   $(\vartheta_1,\vartheta_2) = (1,3)$  & 0.948 & 0.6133  & 0.8376 &  0.6415 & 0.6973 \\
  \hline
   $(\vartheta_1,\vartheta_2) = (1,1)$  & 0.9304  & 0.6079  & 0.4588  &  0.4344 & 0.8775 \\
  \hline
  $(\vartheta_1,\vartheta_2) = (-1,2)$  & 0.8643  & 0.5518  & 0.226 & 0.5238 & 0.5329 \\
  \hline
  $(\vartheta_1,\vartheta_2) = (-3,-2)$  & 0.8454  & 0.61  & 0.8005 & 0.5944  &  0.2997 \\
  \hline
\end{tabular}
\end{small}
\caption{The \ $p$-values of Kolmogorov-Smirnov (KS),  Pearson's chi-squared (PCS), Anderson-Darling (AD) and Jarque-Bera (JB) tests
          for testing normality of the scaled errors.
          The 3rd column contains the \ $p$-values of Kolmogorov-Smirnov test for testing standard normality of the scaled errors.}
\label{tests}
\end{table}

All in all, our simulation results are in accordance with our theoretical ones in Theorem \ref{|lambda_1|>|lambda_2|}.

\section{Proofs}\label{section_proofs}

\noindent
\textbf{Proof of Lemma \ref{LSE1}.}
For each \ $n\in\NN$, \ $(\vartheta_1,\vartheta_2)^\top \in\RR^2$ \ and \ $(x_{-1}, x_0, x_1, \ldots, x_n)^\top \in \RR^{n+2}$, \ we have
 \begin{align*}
  Q_n(x_{-1}, x_0, x_1, \ldots, x_n ; \vartheta_1, \vartheta_2)
  &= \vartheta_1^2 \sum_{k=1}^n x_{k-1}^2 + 2 \vartheta_1 \vartheta_2 \sum_{k=1}^n x_{k-1} x_{k-2}
     + \vartheta_2^2 \sum_{k=1}^n x_{k-2}^2 \\
  &\quad
     - 2 \vartheta_1 \sum_{k=1}^n x_k x_{k-1}
     - 2 \vartheta_2 \sum_{k=1}^n x_k x_{k-2} + \sum_{k=1}^n x_k^2 ,
 \end{align*}
 hence the function \ $\RR^2 \ni (\vartheta_1, \vartheta_2)^\top \mapsto Q_n(x_{-1}, x_0, x_1, \ldots, x_n ; \vartheta_1, \vartheta_2)$ \ is strictly convex if
 \[
   \det\left(\begin{bmatrix}
              \sum_{k=1}^n x_{k-1}^2 & \sum_{k=1}^n x_{k-1} x_{k-2} \\
              \sum_{k=1}^n x_{k-1} x_{k-2} & \sum_{k=1}^n x_{k-2}^2
             \end{bmatrix}\right)
   = \det(G_n(x_{-1}, x_0, x_1, \ldots, x_n))
   > 0 .
 \]
Indeed, if \ $\det(G_n(x_{-1}, x_0, x_1, \ldots, x_n))>0$, \ then \ $ \sum_{k=1}^n x_{k-1}^2 > 0$ \ holds as well, so the matrix
 \[
   G_n(x_{-1}, x_0, x_1, \ldots, x_n)
    = \begin{bmatrix}
              \sum_{k=1}^n x_{k-1}^2 & \sum_{k=1}^n x_{k-1} x_{k-2} \\
              \sum_{k=1}^n x_{k-1} x_{k-2} & \sum_{k=1}^n x_{k-2}^2
       \end{bmatrix}
 \]
 is positive definite yielding that the function
 \ $\RR^2 \ni (\vartheta_1, \vartheta_2)^\top \mapsto Q_n(x_{-1}, x_0, x_1, \ldots, x_n ; \vartheta_1, \vartheta_2)$ \ is strictly convex.
Hence, for each \ $n\in\NN$ \ and for any least squares estimator \ $F_n$ \ of \ $(\vartheta_1, \vartheta_2)^\top$, \ we have if
 \ $(x_{-1}, x_0, x_1, \ldots, x_n)^\top\in D_n$, \ then \ $F_n(x_{-1}, x_0, x_1, \ldots, x_n)$ \ is the
 unique solution of the linear system of equations
 \[
   \begin{cases}
    \vartheta_1 \sum_{k=1}^n x_{k-1}^2 + \vartheta_2 \sum_{k=1}^n x_{k-1} x_{k-2}
    = \sum_{k=1}^n x_k x_{k-1} , \\
    \vartheta_1 \sum_{k=1}^n x_{k-1} x_{k-2} + \vartheta_2 \sum_{k=1}^n x_{k-2}^2
    = \sum_{k=1}^n x_k x_{k-2} .
   \end{cases}
 \]
Consequently, on the set \ $D_n$, \ any least squares estimator \ $F_n$ \ of \ $(\vartheta_1, \vartheta_2)^\top$ \ has the form given in \eqref{LSE_vartheta_1_vartheta_2}.
\proofend

\smallskip
\noindent
\textbf{Proof of Lemma \ref{LEMMA_LSE_exist}.}
It is enough to check that for each \ $n\in\NN$ \ with \ $n\geq 3$, \ we have
 \[
   \PP(\Omega_n) = \PP\left( \det\left( \sum_{k=1}^n  \begin{bmatrix}
                                                        X_{k-1} \\
                                                        X_{k-2} \\
                                                      \end{bmatrix}
                                                       \begin{bmatrix}
                                                        X_{k-1} \\
                                                        X_{k-2} \\
                                                      \end{bmatrix}^\top
                       \right) > 0  \right)=1,
 \]
 since then, by Lemma \ref{LSE1},
 \[
   F_n(X_{-1},X_0,X_1,\ldots,X_n)
        = \begin{bmatrix}
            \hvartheta_1^{(n)} \\
            \hvartheta_2^{(n)} \\
          \end{bmatrix}
    \qquad \text{on the event \ $\Omega_n$,}
  \]
 where \ $F_n : \RR^{n+2}\to\RR^2$ \ is a measurable function satisfying \eqref{LSE_vartheta_1_vartheta_2} on the set \ $D_n$.
\ Here for each \ $n\in\NN$, \ we get
 \begin{align*}
   \det\left( \sum_{k=1}^n \begin{bmatrix}
                               X_{k-1} \\
                               X_{k-2} \\
                           \end{bmatrix}
                           \begin{bmatrix}
                              X_{k-1} \\
                              X_{k-2} \\
                            \end{bmatrix}^\top \right)
    & = \det\left( \begin{bmatrix}
                    \sum_{k=1}^n X_{k-1}^2 & \sum_{k=1}^n X_{k-1}X_{k-2} \\
                    \sum_{k=1}^n X_{k-1}X_{k-2} & \sum_{k=1}^n X_{k-2}^2 \\
                  \end{bmatrix}
            \right) \\
    & = \sum_{k=1}^n X_{k-1}^2 \sum_{k=1}^n X_{k-2}^2
        - \left( \sum_{k=1}^n X_{k-1}X_{k-2} \right)^2.
 \end{align*}
By Cauchy-Schwartz's inequality, we have
 \[
    \sum_{k=1}^n X_{k-1}^2 \sum_{k=1}^n X_{k-2}^2
        - \left( \sum_{k=1}^n X_{k-1}X_{k-2} \right)^2\geq 0, \qquad n\in\NN,
 \]
 and equality holds if and only if the random vectors
 \[
   \begin{bmatrix}
     X_0 & X_1 & \ldots & X_{n-1} \\
   \end{bmatrix}^\top
   \qquad \text{and}\qquad
   \begin{bmatrix}
     X_{-1} & X_0 & \ldots & X_{n-2} \\
   \end{bmatrix}^\top
 \]
 are linearly dependent, i.e., there exist \ $K,L\in\RR$ \ (depending on \ $\omega\in\Omega$) \ such that \ $K^2 + L^2 > 0$ \ and
 \begin{align}\label{help_lin_dep}
    K \begin{bmatrix}
     X_0 & X_1 & \ldots & X_{n-1} \\
   \end{bmatrix}^\top
   + L \begin{bmatrix}
     X_{-1} & X_0 & \ldots & X_{n-2} \\
   \end{bmatrix}^\top
  = \begin{bmatrix}
     0 & 0 & \ldots & 0 \\
   \end{bmatrix}^\top .
 \end{align}
In what follows let us suppose that \ $n\in\NN$ \ and \ $n\geq 3$.
\ If \eqref{help_lin_dep} holds with \ $K=0$ \ (yielding \ $L\ne 0$) \ or \ $L=0$ \ (yielding \ $K\ne 0$), then \ $X_1=0$, \ which can occur only with probability zero,
 since \ $X_1$ \ is absolutely continuous (following from the facts that \ $X_1=\vartheta_1 X_0 + \vartheta_2 X_{-1}+Z_1$ \
 and \ $Z_1$ \ is absolutely continuous).
So if \eqref{help_lin_dep} holds, then, using that the probability space \ $(\Omega,\cF_\infty,\PP)$
 \ is complete, we have \ $\{K=0\}\cup\{L=0\}\in \cF_\infty$ \ and \ $\PP(\{K=0\}\cup\{L=0\})=0$.
\ Further, if \eqref{help_lin_dep} holds with \ $X_{-1}=0$, \ then \ $X_1=0$, \ since if \ $L=0$, \ then \ $K>0$ \ and \ $X_0=X_1=\cdots=X_{n-1}=0$;
 \ and if \ $L>0$ \ and \ $K=0$, \ then \ $X_{-1} = X_0 = X_1 = \cdots = X_{n-2} = 0$; \ and if \ $L>0$ \ and \ $K>0$, \ then
 \ $X_0=0$, \ yielding \ $X_1=0$.
\ Similarly, if \eqref{help_lin_dep} holds with \ $X_0=0$, \ then \ $X_1=0$.
\ So, using again the absolute continuity of \ $X_1$, \ if \eqref{help_lin_dep} holds, then
 \ $\PP(\{X_{-1}=0\} \cup \{X_0=0\})=0$.

\noindent Consequently,
 \begin{align*}
  &\PP(\overline\Omega_n)
   =\PP\left( \det\left( \sum_{k=1}^n \begin{bmatrix}
                                        X_{k-1} \\
                                        X_{k-2} \\
                                      \end{bmatrix}
                                      \begin{bmatrix}
                                        X_{k-1} \\
                                        X_{k-2} \\
                                      \end{bmatrix}^\top
               \right) = 0  \right) \\
  &\leq
    \PP\left(
    \begin{bmatrix}
     X_0 \\
     X_1 \\
     \vdots \\
     X_{n-1} \\
   \end{bmatrix}
   = \left(-\frac{L}{K}\right)
    \begin{bmatrix}
     X_{-1} \\
     X_0 \\
     \vdots \\
     X_{n-2} \\
   \end{bmatrix},
   \;\;\; K\ne 0,\, L\ne 0, \, X_{-1}\ne 0, \, X_0\ne 0
      \right) \\
  & = \PP\left( X_{k-1} = \left( -\frac{L}{K}\right)^k X_{-1}, \, k=1,\ldots,n;\;\; K\ne 0,\, L\ne 0, \, X_{-1}\ne 0, \, X_0\ne 0  \right)\\
  & \leq \PP\left( X_{k-1} = \left( -\frac{L}{K}\right)^k X_{-1}, \, k=1,\ldots,n;\;\; \frac{X_0}{X_{-1}} = -\frac{L}{K},
                            \, K\ne 0, \, L\ne 0, \, X_{-1}\ne 0, \, X_0\ne 0  \right)\\
  & \leq \PP\left( X_{k-1} = \left( \frac{X_0}{X_{-1}}\right)^k X_{-1}, \, k=1,\ldots,n, \, X_{-1}\ne 0, \, X_0\ne 0 \right) \\
  & \leq \PP\left( X_1 = \left(\frac{X_0}{X_{-1}}\right)^2 X_{-1}, \, X_{-1}\ne 0, \, X_0\ne 0\right)\\
  &= \PP\left( \vartheta_1 X_0 + \vartheta_2 X_{-1} + Z_1 = \left(\frac{X_0}{X_{-1}}\right)^2 X_{-1},
                       \, X_{-1}\ne 0, \, X_0\ne 0 \right) \\
  &= \PP\left( Z_1 =  \frac{X_0^2}{X_{-1}} - \vartheta_1 X_0 - \vartheta_2 X_{-1}, \, X_{-1}\ne 0, \, X_0\ne 0 \right)\\
  & = \int_{\RR^2\setminus\{(0,0)\}} \PP\left( Z_1 =  \frac{x_0^2}{x_{-1}} - \vartheta_1 x_0 - \vartheta_2 x_{-1} \right)
             F_{X_{-1}, X_0}(\dd x_{-1},\dd x_0)
    = 0,
 \end{align*}
 since \ $Z_1$ \ and \ $(X_{-1},X_0)$ \ are \ $\PP$-independent, and \ $Z_1$ \ is absolutely continuous, where
 \ $F_{X_{-1}, X_0}$ \ denotes the distribution function of \ $(X_{-1},X_0)$.
\proofend

\smallskip

\noindent
\textbf{Proof of Theorem \ref{|lambda_1|>|lambda_2|}.}
We divide the proof into six steps.

{\sl Step 1 (a decomposition of the left-hand side of \eqref{conv_|lambda_1|>|lambda_2|}):}
For each \ $n \in \NN$ \ with \ $n\geq 3$, \ by \eqref{LSE_diff_formula} and \eqref{langle_bM_rangle_n}, we have
 \begin{align*}
  &\begin{bmatrix}
    \left(\sum_{k=1}^n X_{k-1}^2\right)^{1/2}
     & \frac{\sum_{k=1}^n X_{k-1} X_{k-2}}{\left(\sum_{k=1}^n X_{k-1}^2\right)^{1/2}} \\[3mm]
    \frac{\sum_{k=1}^n X_{k-1} X_{k-2}}{\left(\sum_{k=1}^n X_{k-2}^2\right)^{1/2}}
     & \left(\sum_{k=1}^n X_{k-2}^2\right)^{1/2}
   \end{bmatrix}
   \begin{bmatrix}
    \hvartheta_1^{(n)} - \vartheta_1 \\
    \hvartheta_2^{(n)} - \vartheta_2
   \end{bmatrix} =
 \end{align*}
 \begin{align*}
  &= \begin{bmatrix}
      \left(\sum_{k=1}^n X_{k-1}^2\right)^{-1/2} & 0 \\
      0 & \left(\sum_{k=1}^n X_{k-2}^2\right)^{-1/2}
     \end{bmatrix}
     \begin{bmatrix}
      \sum_{k=1}^n X_{k-1}^2 & \sum_{k=1}^n X_{k-1} X_{k-2} \\
      \sum_{k=1}^n X_{k-1} X_{k-2} & \sum_{k=1}^n X_{k-2}^2
     \end{bmatrix}
     \begin{bmatrix}
      \hvartheta_1^{(n)} - \vartheta_1 \\
      \hvartheta_2^{(n)} - \vartheta_2
     \end{bmatrix} \\
  &= (\sigma^{-1} \bA_n) (\sigma^2 \langle\bM\rangle_n) (\langle\bM\rangle_n^{-1} \bM_n)
   = \sigma \bA_n \bM_n ,
 \end{align*}
 where
 \begin{align}\label{help_An}
   \bA_n := \sigma
            \begin{bmatrix}
             \left(\sum_{k=1}^n X_{k-1}^2\right)^{-1/2} & 0 \\
             0 & \left(\sum_{k=1}^n X_{k-2}^2\right)^{-1/2}
            \end{bmatrix} , \qquad n \in \NN ,
 \end{align}
 on an event having probability one.
Indeed, by Lemma \ref{LEMMA_LSE_exist}, for each \ $n\in\NN$ \ with \ $n\geq 3$, \ the least squares estimator
 \ $(\hvartheta_1^{(n)}, \hvartheta_2^{(n)})^\top$ \ exists uniquely on the event \ $\Omega_n$ \ (given in \eqref{Omega_n})
 having probability one, and, since \ $X_1$ \ is absolutely continuous (following from the facts that
 \ $X_1 = \vartheta_1 X_0 + \vartheta_2 X_{-1} + Z_1$ \ and \ $Z_1$ \ is absolutely continuous), we have
 \ $\PP(\sum_{k=1}^n X_{k-1}^2 > 0)= \PP(\sum_{k=1}^n X_{k-2}^2 > 0) = 1$ \ for each \ $n\in\NN$ \ with \ $n\geq 3$,
 \ yielding that \ $\bA_n$ \ is well-defined for each \ $n\in\NN$ \ with \ $n\geq 3$ \ $\PP$-almost surely.
Hence, by part (c) of Theorem 3.18 in H\"ausler and Luschgy \cite{HauLus}, in order to prove \eqref{conv_|lambda_1|>|lambda_2|},
 it is enough to show that
 \begin{equation}\label{conv_AnMn}
  \bA_n \bM_n
  \to N \begin{bmatrix} 1 \\ \sign(\lambda_1) \end{bmatrix} \qquad
  \text{$\cF_\infty$-mixing as \ $n \to \infty$.}
 \end{equation}
We are going to apply Theorem \ref{MSLTES} with \ $d:=2$, \ $(\bU_n)_{n\in\ZZ_+} := (\bM_n)_{n\in\ZZ_+}$, \ $(\bB_n)_{n\in\NN} := (\bA_n)_{n\in\NN}$,
 \ $(\bQ_n)_{n\in\NN} := (\lambda_1^{-n} \bI_2)_{n\in\NN}$, \ $(\cF_n)_{n\in\ZZ_+} := (\sigma(X_{-1},X_0,X_1,\ldots,X_n))_{n\in\ZZ_+}$, \ and \ $G := \Omega$
 \ (the random matrix \ $\Beta$, \ the matrix \ $\bP$ \ and the probability measure \ $\mu$ \
 appearing in assumptions (i), (iii) and (iv) of Theorem \ref{MSLTES}, respectively, will be
 chosen later on in Step 3).
Note that \ $\bA_n$ \ is invertible for each \ $n\in\NN$ \ with \ $n\geq 3$ \ $\PP$-almost surely.

{\sl Step 2 (asymptotic behaviour of \ $\bA_n$ \ given in \eqref{help_An} as \ $n\to\infty$):}
In order to check the conditions (i)--(iv) of Theorem \ref{MSLTES} with the choices given at the end of Step 1,
 we need the asymptotic behavior of \ $\bA_n$ \ as \ $n \to \infty$.
\ The vectors
 \[
   \begin{bmatrix}
    \lambda_+ \\
    1
   \end{bmatrix} \qquad \text{and} \qquad
   \begin{bmatrix}
    \lambda_- \\
    1
   \end{bmatrix}
 \]
 are right eigenvectors of the matrix \ $\bvartheta$ \ corresponding to the eigenvalues \ $\lambda_+$ \ and \ $\lambda_-$ \ (given in \eqref{help_lambdak}),
  respectively.
Due to our assumption \ $\vert \lambda_1\vert > \vert \lambda_2\vert$ \ and our notations,
 we have \ $\lambda_+ \ne \lambda_-$, \ hence the matrix \ $\bvartheta$ \ can be written in a Jordan canonical form
 \[
   \bvartheta = \begin{bmatrix}
                 \lambda_+ & \lambda_- \\
                 1 & 1
                \end{bmatrix}
                \begin{bmatrix}
                 \lambda_+ & 0 \\
                 0 & \lambda_-
                \end{bmatrix}
                \begin{bmatrix}
                 \lambda_+ & \lambda_- \\
                 1 & 1
                \end{bmatrix}^{-1} .
 \]
Consequently, for each \ $n \in \ZZ_+$, \ we have
 \begin{equation}\label{bvartheta^n}
  \begin{aligned}
   \bvartheta^n
   &= \begin{bmatrix}
       \lambda_+ & \lambda_- \\
       1 & 1
      \end{bmatrix}
      \begin{bmatrix}
       \lambda_+ & 0 \\
       0 & \lambda_-
      \end{bmatrix}^n
      \begin{bmatrix}
       \lambda_+ & \lambda_- \\
       1 & 1
      \end{bmatrix}^{-1} \\
   &= \frac{1}{\lambda_+-\lambda_-}
      \begin{bmatrix}
       \lambda_+ & \lambda_- \\
       1 & 1
      \end{bmatrix}
      \begin{bmatrix}
       \lambda_+^n & 0 \\
       0 & \lambda_-^n
      \end{bmatrix}
      \begin{bmatrix}
       1 & -\lambda_- \\
       -1 & \lambda_+
      \end{bmatrix} \\
   &= \frac{\lambda_+^n}{\lambda_+-\lambda_-}
      \begin{bmatrix}
       \lambda_+ & - \lambda_- \lambda_+ \\
       1 & - \lambda_-
      \end{bmatrix}
      + \frac{\lambda_-^n}{\lambda_+-\lambda_-}
        \begin{bmatrix}
         - \lambda_- & \lambda_+ \lambda_- \\
         - 1 & \lambda_+
        \end{bmatrix} \\
   &= \frac{\lambda_1^n}{\lambda_1-\lambda_2}
      \begin{bmatrix}
       \lambda_1 & - \lambda_1 \lambda_2 \\
       1 & - \lambda_2
      \end{bmatrix}
      + \frac{\lambda_2^n}{\lambda_1-\lambda_2}
        \begin{bmatrix}
         - \lambda_2 & \lambda_1 \lambda_2 \\
         - 1 & \lambda_1
        \end{bmatrix} .
  \end{aligned}
 \end{equation}
For each \ $n \in \NN$, \ by \eqref{X_n_X_n-1} with \ $k=n$, \ we obtain
 \begin{equation}\label{X_n}
  X_n = \begin{bmatrix} 1 \\ 0 \end{bmatrix}^\top
        \begin{bmatrix} X_n \\ X_{n-1} \end{bmatrix}
      = \begin{bmatrix} 1 \\ 0 \end{bmatrix}^\top
        \bvartheta^n
        \begin{bmatrix} X_0 \\ X_{-1} \end{bmatrix}
        + \sum_{j=1}^n
           \begin{bmatrix} 1 \\ 0 \end{bmatrix}^\top
           \bvartheta^{n-j}
           \begin{bmatrix} Z_j \\ 0 \end{bmatrix} .
 \end{equation}
Hence for each \ $n \in \ZZ_+$, \ by \eqref{bvartheta^n}, we have
 \[
   \begin{bmatrix} 1 \\ 0 \end{bmatrix}^\top \bvartheta^n
   = \frac{\lambda_1^{n+1}}{\lambda_1-\lambda_2}
     \begin{bmatrix} 1 \\ - \lambda_2 \end{bmatrix}^\top
     + \frac{\lambda_2^{n+1}}{\lambda_1-\lambda_2}
       \begin{bmatrix} - 1 \\ \lambda_1 \end{bmatrix}^\top .
 \]
Thus, by \eqref{X_n}, for each \ $n \in \NN$, \ we get
 \begin{align*}
  X_n &= \frac{\lambda_1^{n+1}}{\lambda_1-\lambda_2} (X_0 - \lambda_2 X_{-1})
         + \frac{\lambda_2^{n+1}}{\lambda_1-\lambda_2} (-X_0 + \lambda_1 X_{-1}) \\
      &\quad
         + \frac{\lambda_1}{\lambda_1-\lambda_2} \sum_{j=1}^n \lambda_1^{n-j} Z_j
         + \frac{\lambda_2}{\lambda_1-\lambda_2} \sum_{j=1}^n \lambda_2^{n-j} (-Z_j) .
 \end{align*}
Hence we obtain
 \begin{align}\label{help_new1}
  \begin{split}
   \lambda_1^{-n} X_n
   &= \frac{\lambda_1}{\lambda_1-\lambda_2} (X_0 - \lambda_2 X_{-1})
      + \frac{\lambda_2}{\lambda_1-\lambda_2} (-X_0 + \lambda_1 X_{-1}) \left(\frac{\lambda_2}{\lambda_1}\right)^n \\
   &\quad
      + \frac{\lambda_1}{\lambda_1-\lambda_2} \sum_{j=1}^n \lambda_1^{-j} Z_j
      - \frac{\lambda_2}{\lambda_1-\lambda_2} \left(\frac{\lambda_2}{\lambda_1}\right)^n \sum_{j=1}^n \lambda_2^{-j} Z_j \\
   &\asP \frac{\lambda_1}{\lambda_1-\lambda_2} (X_0 - \lambda_2 X_{-1})
        + \frac{\lambda_1}{\lambda_1-\lambda_2} \sum_{j=1}^\infty \lambda_1^{-j} Z_j
    = Y \qquad \text{as \ $n \to \infty$,}
  \end{split}
 \end{align}
 where the \ $\PP$-almost sure absolute convergence of the series \ $\sum_{j=1}^\infty \lambda_1^{-j} Z_j$ \ follows
 by Lemma 8.1 in H\"ausler and Luschgy \cite{HauLus}.
Indeed, since \ $\EE_\PP(\log^+(\vert Z_1\vert))<\infty$, \ $\lambda_1\in\RR$ \ and \ $|\lambda_1| > 1$,
 \ Lemma 8.1 in H\"ausler and Luschgy \cite{HauLus} yields the \ $\PP$-almost sure absolute convergence of
 \ $\sum_{j=1}^\infty \lambda_1^{-j} Z_j$, \ and we check that \ $\bigl(\frac{\lambda_2}{\lambda_1}\bigr)^n \sum_{j=1}^n \lambda_2^{-j} Z_j \asP 0$
 \ as \ $n \to \infty$.
For this it is enough to verify that
 \ $\sum_{n=1}^\infty \PP\Big(\Big|\Big(\frac{\lambda_2}{\lambda_1}\Big)^n \sum_{j=1}^n \lambda_2^{-j} Z_j\Big| > \vare\Big)<\infty$ \
 for each \ $\vare\in\RR_{++}$.
\ Using that \ $Z_j$, \ $j\in\NN$, \ are i.i.d.\ random variables having zero mean and variance \ $\sigma^2$,
 \ for each \ $\vare \in \RR_{++}$, \ we have
 \begin{align*}
  \sum_{n=1}^\infty
   \PP\left(\left|\left(\frac{\lambda_2}{\lambda_1}\right)^n \sum_{j=1}^n \lambda_2^{-j} Z_j\right| > \vare\right)
  &\leq \frac{1}{\vare^2}
        \sum_{n=1}^\infty
         \EE_\PP\left(\left|\left(\frac{\lambda_2}{\lambda_1}\right)^n \sum_{j=1}^n \lambda_2^{-j} Z_j\right|^2\right) \\
  &= \frac{1}{\vare^2}
     \sum_{n=1}^\infty
      \left(\frac{\lambda_2}{\lambda_1}\right)^{2n}
      \sum_{j=1}^n
       \lambda_2^{-2j} \sigma^2
   < \infty ,
 \end{align*}
 since
 \[
   \sum_{j=1}^n \lambda_2^{-2j}
       = \begin{cases}
            n & \text{if \ $\vert\lambda_2\vert=1$,}\\
            \lambda_2^{-2} \frac{\lambda_2^{-2n} - 1}{\lambda_2^{-2} - 1}
              & \text{if \ $\vert\lambda_2\vert\ne 1$,}
         \end{cases}
 \]
 and if \ $\vert\lambda_2\vert=1$, \ then
 \[
    \limsup_{n\to\infty} \left\vert \frac{ \left(\frac{\lambda_2}{\lambda_1}\right)^{2(n+1)} (n+1) }
                                         {  \left(\frac{\lambda_2}{\lambda_1}\right)^{2n} n } \right\vert
     = \left\vert \frac{\lambda_2}{\lambda_1} \right\vert^2
     = \vert \lambda_1\vert^{-2} < 1,
 \]
 and if \ $\vert\lambda_2\vert\ne 1$, \ then \ $|\lambda_1| > |\lambda_2|$ \ implies that
 \begin{align*}
    &\limsup_{n\to\infty} \left\vert \frac{ \left(\frac{\lambda_2}{\lambda_1}\right)^{2(n+1)} \frac{\lambda_2^{-2(n+1)} - 1}{\lambda_2^{-2} - 1}  }
                                         {  \left(\frac{\lambda_2}{\lambda_1}\right)^{2n}  \frac{\lambda_2^{-2n} - 1}{\lambda_2^{-2} - 1} } \right\vert
      = \left\vert \frac{\lambda_2}{\lambda_1} \right\vert^2
       \limsup_{n\to\infty} \left\vert  \frac{\lambda_2^{-2(n+1)} - 1}{\lambda_2^{-2n} - 1 } \right\vert\\
    &\qquad = \begin{cases}
         \left\vert \frac{\lambda_2}{\lambda_1} \right\vert^2 < 1 & \text{if \ $\vert \lambda_2\vert> 1$,}\\
         \left\vert \frac{\lambda_2}{\lambda_1} \right\vert^2
          \limsup_{n\to\infty} \left\vert  \frac{\lambda_2^{-2} - \lambda_2^{2n}}{1- \lambda_2^{2n} } \right\vert
          = \left\vert \frac{\lambda_2}{\lambda_1} \right\vert^2  \vert \lambda_2\vert^{-2} = \vert \lambda_1\vert^{-2}<1
             & \text{if \ $\vert \lambda_2\vert< 1$.}
       \end{cases}
 \end{align*}
By D'Alambert's criteria, this yields that the series \ $\sum_{n=1}^\infty \big(\frac{\lambda_2}{\lambda_1}\big)^{2n} \sum_{j=1}^n \lambda_2^{-2j}$ \ is convergent.
We note that if, in addition, \ $(X_k)_{k\geq -1}$ \ is purely explosive as well (i.e., \ $\lambda_1,\lambda_2\in\RR$ \ and
 \ $\vert \lambda_1\vert> \vert \lambda_2\vert>1$), \ then \ $\bigl(\frac{\lambda_2}{\lambda_1}\bigr)^n \sum_{j=1}^n \lambda_2^{-j} Z_j \asP 0$
 \ as \ $n \to \infty$ \ follows more easily, since in this case, by Lemma 8.1 in H\"ausler and Luschgy \cite{HauLus}, \ $\sum_{j=1}^\infty \lambda_2^{-j} Z_j$ \
 is absolutely convergent \ $\PP$-a.s., and, since \ $\vert\lambda_2\vert<\vert\lambda_1\vert$, \ we have
 \ $\bigl(\frac{\lambda_2}{\lambda_1}\bigr)^n \sum_{j=1}^n \lambda_2^{-j} Z_j \asP 0 \cdot\sum_{j=1}^\infty \lambda_2^{-j} Z_j=0$ \ as \ $n\to\infty$.

The random variable \ $Y$ \ given in \eqref{conv_X_n} is absolutely continuous.
Indeed,
 \[
   Y =  \frac{\lambda_1}{\lambda_1-\lambda_2} (X_0 - \lambda_2 X_{-1} + \lambda_1^{-1} Z_1)
        + \frac{\lambda_1}{\lambda_1-\lambda_2} \sum_{j=2}^\infty \lambda_1^{-j} Z_j,
 \]
 where the absolute continuity of \ $Z_1$ \ and the independence of \ $Z_1$ \ and \ $(X_{-1},X_0)$ \ yield the absolute continuity
 of \ $X_0 - \lambda_2 X_{-1} + \lambda_1^{-1} Z_1$.
\ Hence, using that \ $\sum_{j=2}^\infty \lambda_1^{-j} Z_j$ \ and \ $X_0 - \lambda_2 X_{-1} + \lambda_1^{-1} Z_1$ \ are independent,
 we have the absolute continuity of \ $Y$.

By \eqref{help_new1}, we get \ $\lambda_1^{-2n} X_n^2 \asP Y^2$ \ as \ $n \to \infty$, \
 and, since \ $\sum_{k=1}^\infty \lambda_1^{2(k-1)} = \infty$ \ (due to \ $\vert \lambda_1\vert > 1$, \ $\lambda_1\in\RR$),
 \ applying the Toeplitz lemma (see, e.g., H\"ausler and Luschgy \cite[Lemma 6.28]{HauLus}), we get
 \[
   \frac{\sum_{k=1}^n X_{k-1}^2}{\sum_{k=1}^n \lambda_1^{2(k-1)}}
     \asP Y^2 \qquad \text{as \ $n \to \infty$.}
 \]
Consequently, using \ $\sum_{k=1}^n \lambda_1^{2(k-1)} = \frac{\lambda_1^{2n}-1}{\lambda_1^2-1}$, \ we
conclude
 \begin{equation}\label{sum_{k=1}^n X_{k-1}^2-}
  \lambda_1^{-2n} \sum_{k=1}^n X_{k-1}^2
  = \lambda_1^{-2n} \frac{\lambda_1^{2n}-1}{\lambda_1^2-1}
    \frac{\sum_{k=1}^n X_{k-1}^2}{\sum_{k=1}^n \lambda_1^{2(k-1)}}
  \asP \frac{1}{\lambda_1^2-1} Y^2 \qquad \text{as \ $n \to \infty$.}
 \end{equation}
In a similar way, we have
 \begin{align}\label{sum_{k=1}^n X_{k-2}^2-}
  \lambda_1^{-2n} \sum_{k=1}^n X_{k-2}^2
  \asP \frac{1}{\lambda_1^2(\lambda_1^2-1)} Y^2 \qquad \text{as \ $n \to \infty$.}
 \end{align}
Indeed, by \eqref{help_new1} and \eqref{sum_{k=1}^n X_{k-1}^2-}, \ we have
 \begin{align*}
   \lambda_1^{-2n} \sum_{k=1}^n X_{k-2}^2
    & = \lambda_1^{-2n} (X_{-1}^2 - X_{n-1}^2) + \lambda_1^{-2n} \sum_{k=1}^n X_{k-1}^2 \\
    & = \lambda_1^{-2n} X_{-1}^2 - \lambda_1^{-2} \lambda_1^{-2(n-1)} X_{n-1}^2
        + \lambda_1^{-2n} \sum_{k=1}^n X_{k-1}^2 \\
    & \asP 0 - \lambda_1^{-2} Y^2 + (\lambda_1^2 - 1)^{-1} Y^2
         =  \frac{1}{\lambda_1^2(\lambda_1^2-1)} Y^2 \qquad \text{as \ $n \to \infty$,}
 \end{align*}
 as desired.

The absolute continuity of \ $Y$ \ implies \ $\PP(Y = 0) = 0$, \ hence, by \eqref{help_An}, \eqref{sum_{k=1}^n X_{k-1}^2-} and \eqref{sum_{k=1}^n X_{k-2}^2-}, we obtain
 \begin{equation}\label{bA_n}
  \begin{aligned}
   \lambda_1^n \bA_n
   &= \sigma
      \begin{bmatrix}
       \left(\lambda_1^{-2n} \sum_{k=1}^n X_{k-1}^2\right)^{-1/2} & 0 \\
       0 & \left(\lambda_1^{-2n} \sum_{k=1}^n X_{k-2}^2\right)^{-1/2}
      \end{bmatrix} \\
   &\asP \sigma
        \begin{bmatrix}
         \Bigl(\frac{Y^2}{\lambda_1^2-1}\Bigr)^{-1/2} & 0 \\
         0 & \Bigl(\frac{Y^2}{(\lambda_1^2-1)\lambda_1^2}\Bigr)^{-1/2}
        \end{bmatrix}
    = \frac{\sigma\sqrt{\lambda_1^2-1}}{|Y|}
      \begin{bmatrix}
       1 & 0 \\
       0 & |\lambda_1|
      \end{bmatrix}
    \qquad \text{as \ $n \to \infty$.}
  \end{aligned}
 \end{equation}

{\sl Step 3 (checking conditions (i) and (iii) of Theorem \ref{MSLTES}):}
Recall that at the end of Step 1 we gave our choices for \ $(\bU_n)_{n\in\ZZ_+}$, \ $(\bB_n)_{n\in\NN}$, \ $(\bQ_n)_{n\in\NN}$, \ $(\cF_n)_{n\in\ZZ_+}$ \
 and \ $G$ \ in Theorem \ref{MSLTES}, the random matrix \ $\Beta$ \ and the matrix \ $\bP$ \ appearing in assumptions (i) and (iii) of Theorem \ref{MSLTES},
 respectively, will be chosen below.
Applying \eqref{bA_n}, we obtain
 \begin{equation}\label{bA_n^{-1}}
  \begin{aligned}
   \lambda_1^{-n} \bA_n^{-1} = (\lambda_1^n \bA_n)^{-1}
   &\asP \left(\frac{\sigma\sqrt{\lambda_1^2-1}}{|Y|}
        \begin{bmatrix}
         1 & 0 \\
         0 & |\lambda_1|
        \end{bmatrix}\right)^{-1} \\
   &\phantom{\asP\;}
     = \frac{|Y|}{\sigma\sqrt{\lambda_1^2-1}}
      \begin{bmatrix}
       1 & 0 \\
       0 & |\lambda_1|^{-1}
       \end{bmatrix}
      \qquad \text{as \ $n \to \infty$.}
  \end{aligned}
 \end{equation}
Hence, since almost sure convergence yields convergence in probability,
 we obtain that condition (i) of Theorem \ref{MSLTES} holds with
 \[
   \Beta := \frac{|Y|}{\sigma\sqrt{\lambda_1^2-1}}
            \begin{bmatrix}
             1 & 0 \\
             0 & |\lambda_1|^{-1}
            \end{bmatrix} ,
 \]
 which is invertible if and only if \ $Y \ne 0$, \ so \ $\PP(\exists \, \Beta^{-1}) = \PP(Y \ne 0) = 1$.
\ Here \ $\Beta$ \ is \ $\cF_\infty$-measurable, since \ $Y$ \ is \ $\cF_\infty$-mesaurable.
Indeed, the series \ $\sum_{j=1}^\infty \lambda_1^{-j} Z_j$ \ converges \ $\PP$-a.s. (see Step 2),
 \ $\sum_{j=1}^n \lambda_1^{-j} Z_j$ \ is \ $\cF_\infty$-measurable for all \ $n\in\NN$ \
 (due to \ $\sigma(X_{-1},X_0,X_1,\ldots,X_n) = \sigma(X_{-1},X_0,Z_1,\ldots,Z_n)$, \ $n\in\NN$),
 and the probability space \ $(\Omega,\cF_\infty,\PP)$ \ is complete (see, e.g., Cohn \cite[Corollary 2.2.3]{Coh}).
Moreover, by \eqref{bA_n^{-1}}, for every \ $r \in \NN$, \ we have
 \[
   \bA_n \bA_{n-r}^{-1} = \lambda_1^{-r} (\lambda_1^n \bA_n) (\lambda_1^{-(n-r)} \bA_{n-r}^{-1})
   \asP \lambda_1^{-r} \Beta^{-1} \Beta = \lambda_1^{-r} \bI_2 \qquad \text{as \ $n \to \infty$,}
 \]
 hence we obtain that condition (iii) of Theorem \ref{MSLTES} holds with \ $\bP := \lambda_1^{-1} \bI_2$
 \ (of which the spectral radius is strictly less than 1, since \ $\vert \lambda_1\vert>1$).

{\sl Step 4 (checking condition (ii) of Theorem \ref{MSLTES}):}
Condition (ii) of Theorem \ref{MSLTES} with the earlier given choices (see the end of Step 1) holds if and only if
 \ $\bigl(\lambda_1^{-n} M_n^{(j)}\bigr)_{n\in\ZZ_+}$ \ is stochastically bounded
  for each \ $j \in \{1, 2\}$, \ where \ $(M_n^{(j)})_{n\in\ZZ_+} := (\be_j^\top \bM_n)_{n\in\ZZ_+}$ \ for \ $j \in \{1, 2\}$.
\ Indeed, \ $\PP(G\cap\{\exists \; \Beta^{-1}\}) = \PP(Y\ne 0)=1$ \ yielding that \ $\PP_{G\cap\{\exists \; \Beta^{-1}\}} = \PP$, \ and
 for each \ $n\in\ZZ_+$ \ and \ $K\in(0,\infty)$, \ we have
 \begin{align*}
  \PP(\Vert \lambda_1^{-n} \bM_n\Vert > K)
     \leq  \PP\left(\vert \lambda_1^{-n} M^{(1)}_n\vert  > \frac{K}{\sqrt{2}}\right)
           + \PP\left(\vert \lambda_1^{-n} M^{(2)}_n\vert > \frac{K}{\sqrt{2}}\right),
 \end{align*}
 and
 \[
   \PP( \vert \lambda_1^{-n} M_n^{(j)}\vert > K) \leq \PP(\Vert \lambda_1^{-n} \bM_n\Vert > K),
     \qquad j=1,2.
 \]
By \eqref{langle_bM_rangle_n}, the process \ $(M_n^{(1)})_{n\in\ZZ_+}$ \ is a square integrable martingale with respect to the filtration
 \ $(\cF_n)_{n\in\ZZ_+}$ \ and it has a  quadratic characteristic process
 \[
   \langle M^{(1)} \rangle_n
   = \sigma^{-2} \sum_{k=1}^n X_{k-1}^2 , \qquad n \in \NN ,
 \]
 with \ $\langle M^{(1)}\rangle_0=0$.
\ For each \ $n \in \NN$ \ and \ $K \in (0, \infty)$, \ by Lenglart's inequality (see Corollary \ref{Lenglart_Y}), we get
 \[
   \PP(|\lambda_1^{-n} M_n^{(1)}| \geq K)
   = \PP(|M_n^{(1)}|^2 \geq K^2 \lambda_1^{2n})
   \leq \frac{1}{K}
        + \PP(\langle M^{(1)} \rangle_n > K \lambda_1^{2n}) ,
 \]
 so that for each \ $K \in (0, \infty)$, \ we have
 \begin{align}\label{help_Lenglart}
   \sup_{n\in\NN} \PP(|\lambda_1^{-n} M_n^{(1)}| \geq K)
   \leq \frac{1}{K}
        + \sup_{n\in\NN} \PP(\lambda_1^{-2n} \langle M^{(1)} \rangle_n > K) .
 \end{align}
By \eqref{sum_{k=1}^n X_{k-1}^2-}, we get
 \[
   \lambda_1^{-2n} \langle M^{(1)} \rangle_n
   = \lambda_1^{-2n} \sigma^{-2} \sum_{k=1}^n X_{k-1}^2
   \asP \frac{Y^2}{(\lambda_1^2-1)\sigma^2}
   \qquad \text{as \ $n \to \infty$,}
 \]
 hence \ $\bigl(\lambda_1^{-2n} \langle M^{(1)} \rangle_n\bigr)_{n\in\ZZ_+}$ \ is stochastically bounded, i.e.,
 \[
   \lim_{K\to\infty} \sup_{n\in\NN} \PP(\vert \lambda_1^{-2n} \langle M^{(1)} \rangle_n \vert > K) = 0.
 \]
Consequently, by \eqref{help_Lenglart}, \ $\lim_{K\to\infty} \sup_{n\in\NN} \PP(|\lambda_1^{-n} M_n^{(1)}| \geq K) = 0$, \ i.e.,
 \ $\bigl(\lambda_1^{-n} M_n^{(1)}\bigr)_{n\in\ZZ_+}$ \ is stochastically bounded.
Using \eqref{sum_{k=1}^n X_{k-2}^2-}, in a similar way, one can check that the process
 \ $\bigl(\lambda_1^{-n} M_n^{(2)}\bigr)_{n\in\ZZ_+}$ \ is stochastically bounded,
 and we conclude that condition (ii) of Theorem \ref{MSLTES} holds.

{\sl Step 5 (checking condition (iv) of Theorem \ref{MSLTES}):}
In order to check condition (iv) of Theorem \ref{MSLTES}, let us observe that the square integrable martingale \ $(\bM_n)_{n\in\ZZ_+}$ \ has conditional Gaussian increments with respect to the filtration \ $(\cF_n)_{n\in\ZZ_+}$, \ since for each \ $n \in \NN$, \ the conditional distribution of \ $\Delta\bM_n = \bM_n - \bM_{n-1} = \frac{1}{\sigma^2} Z_n \begin{bmatrix} X_{n-1} \\ X_{n-2} \end{bmatrix}$ \ given \ $\cF_{n-1}$ \ is
 \[
   \cN_2\left(\bzero, \frac{1}{\sigma^4}
               \EE_\PP(Z_n^2) \begin{bmatrix} X_{n-1} \\ X_{n-2} \end{bmatrix}
               \begin{bmatrix} X_{n-1} \\ X_{n-2} \end{bmatrix}^\top\right)
   = \cN_2\left(\bzero, \frac{1}{\sigma^2}
                 \begin{bmatrix} X_{n-1} \\ X_{n-2} \end{bmatrix}
                 \begin{bmatrix} X_{n-1} \\ X_{n-2} \end{bmatrix}^\top\right)
   = \cN_2(\bzero, \Delta \langle\bM\rangle_n) ,
 \]
 where the last equality follows by \eqref{langle_bM_rangle_n}.
More precisely, using the notations and results of H\"ausler and Luschgy \cite{HauLus}, for each \ $n\in\NN$, \ the conditional
 distribution \ $\PP^{\Delta \bM_n \mid \cF_{n-1}}$ \ of \ $\Delta\bM_n$ \ given \ $\cF_{n-1}$ \ can be calculated as follows:
 \begin{align*}
  \PP^{\Delta \bM_n \mid \cF_{n-1}}
   & = \PP^{\sigma^{-2} Z_n [X_{n-1} ,\, X_{n-2} ]^\top \mid \cF_{n-1}}
     = \PP^{ g(X_{n-2},X_{n-1},Z_n) \mid \cF_{n-1}}
     = \left(  \PP^{ (X_{n-2},X_{n-1},Z_n) \mid \cF_{n-1}} \right)^g \\
   & = \left( \delta_{(X_{n-2},X_{n-1})} \otimes \PP^{ Z_n \mid \cF_{n-1}}  \right)^g
     = \left( \delta_{(X_{n-2},X_{n-1})}  \otimes \PP^{Z_n} \right)^g,
 \end{align*}
 where \ $g:\RR^3 \to \RR^2$, \ $g(x_1,x_2,z):= \sigma^{-2} z [x_2,\, x_1]^\top$, \ $(x_1,x_2,z)\in\RR^3$, \
 \ $\PP^{Z_n}$ \  denotes the distribution of \ $Z_n$ \ under \ $\PP$, \ $\delta_{(X_{n-2},X_{n-1})}$ \ is the Dirac Markov kernel
 corresponding to \ $(X_{n-2},X_{n-1})$, \ and we used part (a) of Lemma A.5 in H\"ausler and Luschgy \cite{HauLus},
 parts (b) and \ $(c)$ \ of Lemma A.4 in H\"ausler and Luschgy \cite{HauLus},
 the independence of \ $Z_n$ \ and \ $\cF_{n-1}$, \ and the \ $\cF_{n-1}$-measurability of \ $X_{n-1}$ \ and \ $X_{n-2}$.
\ Hence for each \ $n\in\NN$, \ $\omega\in\Omega$, \ and \ $B\in\cB(\RR^2)$, \ we have
 \begin{align*}
  \PP^{\Delta \bM_n \mid \cF_{n-1}}(\omega, B)
     & = (\delta_{(X_{n-2},X_{n-1})} \otimes \PP^{Z_n}) \big(\omega, \{ (x_1,x_2,z)^\top \in\RR^3 : g(x_1,x_2,z)\in B \}\big) \\
     & = \PP^{Z_n} \big( \{z\in\RR : \sigma^{-2}z [X_{n-1}(\omega) ,\, X_{n-2}(\omega) ]^\top \in B  \}\big) \\
     &  = \PP\big( \{ \widetilde\omega\in\Omega : \sigma^{-2}Z_n(\widetilde\omega) [X_{n-1}(\omega) ,\, X_{n-2}(\omega) ]^\top \in B  \}\big)\\
     &= \PP^{\bzeta_{n,\omega}}(B),
 \end{align*}
 where \ $\bzeta_{n,\omega}$ \ is an \ $\RR^2$-valued random variable having distribution
 \[
   \cN_2\left(\bzero, \sigma^{-2}
                 \begin{bmatrix} X_{n-1}(\omega) \\ X_{n-2}(\omega) \end{bmatrix}
                  \begin{bmatrix} X_{n-1}(\omega) \\ X_{n-2}(\omega) \end{bmatrix}^\top\right)
     = \cN_2(\bzero, \Delta \langle\bM\rangle_n(\omega)),
 \]
 as desired.
Consequently, for each \ $\btheta \in \RR^2$ \ and \ $n \in \NN$, \ we obtain
 \begin{align*}
  &\EE_\PP\biggl(\exp\bigl\{\ii \left\langle\btheta, \bA_n \Delta\bM_n\right\rangle\bigr\} \,\bigg|\, \cF_{n-1}\biggr)
    = \EE_\PP\biggl(\exp\biggl\{\ii \left\langle \bA_n^\top \btheta, \Delta\bM_n  \right\rangle\biggr\} \,\bigg|\, \cF_{n-1}\biggr)\\
  & = \exp\biggl\{- \frac{1}{2} \Big\langle (\Delta \langle\bM\rangle_n) \bA_n^\top \btheta, \bA_n^\top \btheta \Big\rangle \biggr\}
    = \exp\biggl\{- \frac{1}{2} \btheta^\top \bA_n  (\Delta \langle\bM\rangle_n) \bA_n \btheta\biggr\} ,
 \end{align*}
 where, at the last equality, we used that \ $\bA_n$ \ is symmetric.
Using \eqref{langle_bM_rangle_n} and \eqref{help_new1}, we get
 \begin{equation}\label{Delta<bM>_n}
  \begin{aligned}
   \lambda_1^{-2n}
   \Delta \langle\bM\rangle_n
   &= \sigma^{-2} \lambda_1^{-2n}
      \begin{bmatrix} X_{n-1} \\ X_{n-2} \end{bmatrix}
      \begin{bmatrix} X_{n-1} \\ X_{n-2} \end{bmatrix}^\top
    = \sigma^{-2} \lambda_1^{-2n}
      \begin{bmatrix} X_{n-1}^2 & X_{n-1} X_{n-2} \\ X_{n-1} X_{n-2} & X_{n-2}^2 \end{bmatrix} \\
   &\asP \sigma^{-2} Y^2
        \begin{bmatrix}
         \lambda_1^{-2} & \lambda_1^{-3} \\
         \lambda_1^{-3} & \lambda_1^{-4}
        \end{bmatrix}
    \qquad \text{as \ $n \to \infty$.}
  \end{aligned}
 \end{equation}
Applying \eqref{bA_n} and \eqref{Delta<bM>_n}, since \ $\PP(Y=0)=0$, \ we get
 \begin{align*}
  \bA_n  (\Delta \langle\bM\rangle_n) \bA_n
  &= (\lambda_1^n \bA_n)  (\lambda_1^{-2n} \Delta \langle\bM\rangle_n) (\lambda_1^n \bA_n) \\
  &\asP \frac{(\lambda_1^2-1)\sigma^2}{Y^2} \cdot \frac{Y^2}{\sigma^2}
       \begin{bmatrix}
        1 & 0 \\
        0 & |\lambda_1|
       \end{bmatrix}
       \begin{bmatrix}
        \lambda_1^{-2} & \lambda_1^{-3} \\
        \lambda_1^{-3} & \lambda_1^{-4}
       \end{bmatrix}
       \begin{bmatrix}
        1 & 0 \\
        0 & |\lambda_1|
       \end{bmatrix} \\
  &\phantom{\asP}= \frac{\lambda_1^2-1}{\lambda_1^2}
     \begin{bmatrix}
      1 & \sign(\lambda_1) \\
      \sign(\lambda_1) & 1
     \end{bmatrix}
    \qquad \text{as \ $n \to \infty$.}
 \end{align*}
Hence, for each \ $\btheta \in \RR^2$, \ we have
 \begin{align*}
  \EE_\PP\biggl(\exp\bigl\{\ii \left\langle\btheta, \bA_n \Delta\bM_n\right\rangle\bigr\} \,\bigg|\, \cF_{n-1}\biggr)
  &\asP \exp\left\{-\frac{\lambda_1^2-1}{2\lambda_1^2} \btheta^\top
                   \begin{bmatrix}
                    1 & \sign(\lambda_1) \\
                    \sign(\lambda_1) & 1
                   \end{bmatrix} \btheta\right\} \\
  &\phantom{\asP} = \int_{\RR^2} \ee^{\ii \langle\btheta, \bx\rangle} \,
             \PP^{\bzeta}(\dd\bx)
 \end{align*}
 as \ $n \to \infty$, \ where \ $\bzeta$ \ is an \ $\RR^2$-valued random variable having distribution
 \[
  \cN_2\left(\bzero, \frac{\lambda_1^2-1}{\lambda_1^2}
                   \begin{bmatrix}
                    1 & \sign(\lambda_1) \\
                    \sign(\lambda_1) & 1
                   \end{bmatrix}\right).
 \]
Hence we obtain that condition (iv) of Theorem \ref{MSLTES} holds with \ $\mu := \PP^{\bzeta}$, \ since
 \begin{align*}
  \int_{\RR^2} \log^+(\Vert \bx\Vert) \,\mu(\dd\bx)
    & = \int_{\{ \bx\in \RR^2 : \Vert \bx\Vert\geq 1\}} \log(\Vert \bx\Vert) \,\mu(\dd\bx)
      \leq \int_{\{\bx\in \RR^2 : \Vert \bx\Vert\geq 1\}} \Vert \bx\Vert \,\mu(\dd\bx) \\
    & \leq \int_{\{\bx=(x_1,x_2)^\top\in \RR^2 : \Vert \bx\Vert\geq 1\}} (\vert x_1\vert + \vert x_2\vert) \,\mu(\dd\bx)
      <\infty
 \end{align*}
 due to the fact that all the mixed moments of \ $\mu$ \ (being a \ $2$-dimensional normal distribution) are finite.

{\sl Step 6 (application of Theorem \ref{MSLTES}):}
Using Steps 1--5, we can apply Theorem \ref{MSLTES} with our choices given at the end of Step 1 and
 in Steps 3 and 5 (for \ $\Beta$, \ $\bP$ \ and \ $\mu$), \ and we obtain
 \begin{equation}\label{conv_AnMn_MSLTES}
  \bA_n \bM_n
  \to \sum_{j=0}^\infty \begin{bmatrix}
                         \lambda_1^{-j} & 0 \\
                         0 & \lambda_1^{-j}
                        \end{bmatrix}
                        \bN_j \qquad
  \text{$\cF_\infty$-mixing as \ $n \to \infty$,}
 \end{equation}
 where \ $(\bN_j)_{j\in\ZZ_+}$ \ is a sequence of \ $\PP$-independent and identically distributed \ $\RR^2$-valued random vectors being
 \ $\PP$-independent of \ $\cF_\infty$ \ with \ $\PP(\bN_0 \in B) = \mu(B)$ \ for all \ $B\in\cB(\RR^2)$, \ and the series in \eqref{conv_AnMn_MSLTES}
 converges \ $\PP$-almost surely.
The distribution of the limit random variable in \eqref{conv_AnMn_MSLTES} can be written in the form
 \[
   \sum_{j=0}^\infty \begin{bmatrix}
                      \lambda_1^{-j} & 0 \\
                      0 & \lambda_1^{-j}
                     \end{bmatrix}
                     \bN_j
   = \sum_{j=0}^\infty \lambda_1^{-j} \bN_j
   \distre \begin{bmatrix} N \\ \sign(\lambda_1) N \end{bmatrix}
   = N \begin{bmatrix} 1 \\ \sign(\lambda_1) \end{bmatrix} ,
 \]
 where \ $N$ \ is \ $\PP$-independent of \ $\cF_\infty$ \ and \ $N \distre \cN(0,1)$, \ since the \ $2$-dimensional random vector
 \ $\sum_{j=0}^\infty \lambda_1^{-j}\bN_j$ \ has a \ $2$-dimensional normal distribution with covariance matrix
 \[
   \left(\sum_{j=0}^\infty \lambda_1^{-2j}\right)
    \frac{\lambda_1^2-1}{\lambda_1^2}
         \begin{bmatrix}
          1 & \sign(\lambda_1) \\
          \sign(\lambda_1) & 1
         \end{bmatrix}
   = \begin{bmatrix}\
      1 & \sign(\lambda_1) \\
      \sign(\lambda_1) & 1
     \end{bmatrix} .
 \]
Consequently, we conclude \eqref{conv_AnMn}, hence, as it was explained,
 the convergence \eqref{conv_|lambda_1|>|lambda_2|} follows from part (c) of Theorem 3.18 in H\"ausler and Luschgy \cite{HauLus}.
\proofend

\smallskip

\noindent{\bf Proof of Corollary \ref{Cor1}.}
In the proof of Theorem \ref{|lambda_1|>|lambda_2|}, we showed that Theorem \ref{MSLTES} can be applied with the choices
 \ $d:=2$, \ $(\bU_n)_{n\in\ZZ_+}:=(\bM_n)_{n\in\ZZ_+}$, \ $(\bB_n)_{n\in\NN} := (\bA_n)_{n\in\NN}$, \ where \ $\bA_n$, $n\in\NN$, \ is given in \eqref{help_An},
 \ $(\bQ_n)_{n\in\NN} := (\lambda_1^{-n}\bI_2)_{n\in\NN}$, \ $(\cF_n)_{n\in\ZZ_+} := (\sigma(X_{-1},X_0,X_1,\ldots,X_n))_{n\in\ZZ_+}$, \ $G:=\Omega$, \
 and \ $\Beta$ \ and \ $\bP$ \ are given in Step 3 in the proof of Theorem \ref{|lambda_1|>|lambda_2|}.
So, by \eqref{conv_QU}, we have the statement.
\proofend

\smallskip

\noindent{\bf Proof of Corollary \ref{Cor_detMn}.}
First, we prove
 \begin{align}\label{sum_{k=1}^n X_{k-1} X_{k-2}-}
  \lambda_1^{-2n} \sum_{k=1}^n X_{k-1} X_{k-2}
  \asP \frac{1}{(\lambda_1^2 - 1)\lambda_1} Y^2 \qquad \text{as \ $n \to \infty$.}
 \end{align}
By \eqref{help_new1}, we have
 \[
    \lambda_1^{-n} X_n \asP Y \qquad \text{as \ $n\to\infty$,}\qquad \text{and} \qquad \lambda_1^{-(n-1)} X_{n-1} \asP Y \qquad \text{as \ $n\to\infty$,}
 \]
 so
 \[
    \lambda_1^{-n} (X_n + X_{n-1}) = \lambda_1^{-n}  X_n + \lambda_1^{-1} \lambda_1^{-(n-1)} X_{n-1}
       \asP  \left( 1+ \lambda_1^{-1}\right) Y
        \qquad \text{as \ $n\to\infty$.}
 \]
Since \ $\sum_{k=1}^\infty \lambda_1^{2(k-1)} = \infty$ \ (due to \ $\vert \lambda_1\vert > 1$, \ $\lambda_1\in\RR$),
 \ applying the Toeplitz lemma (see, e.g., H\"ausler and Luschgy \cite[Lemma 6.28]{HauLus}), we get
 \[
   \frac{\sum_{k=1}^n (X_{k-1} + X_{k-2})^2 }
        {\sum_{k=1}^n \lambda_1^{2(k-1)} }
        \asP (1+\lambda_1^{-1})^2 Y^2
       \qquad \text{as \ $n\to\infty$.}
 \]
Since \ $\sum_{k=1}^n \lambda_1^{2(k-1)} = \frac{\lambda_1^{2n} -1}{\lambda_1^2 -1}$, \ we have
 \[
    \lambda_1^{-2n} \sum_{k=1}^n (X_{k-1} + X_{k-2})^2
      \asP \frac{(\lambda_1^{-1} + 1)^2}{\lambda_1^2 -1} Y^2
       \qquad \text{as \ $n\to\infty$,}
 \]
 and, by \eqref{sum_{k=1}^n X_{k-1}^2-} and \eqref{sum_{k=1}^n X_{k-2}^2-},
 \begin{align*}
   \lambda_1^{-2n} \sum_{k=1}^n X_{k-1}X_{k-2}
      &=  \frac{1}{2}\lambda_1^{-2n}\left[  \sum_{k=1}^n (X_{k-1} + X_{k-2})^2 -  \sum_{k=1}^n X_{k-1}^2 -  \sum_{k=1}^n X_{k-2}^2   \right] \\
      &\asP \frac{1}{2}\left[  \frac{(\lambda_1^{-1} + 1)^2}{\lambda_1^2 -1} - \frac{1}{\lambda_1^2 -1} -  \frac{1}{(\lambda_1^2 -1)\lambda_1^2} \right] Y^2
           = \frac{1}{(\lambda_1^2 -1)\lambda_1}Y^2
 \end{align*}
 as \ $n\to\infty$, \ yielding \eqref{sum_{k=1}^n X_{k-1} X_{k-2}-}.

Finally, \eqref{langle_bM_rangle_n}, \eqref{sum_{k=1}^n X_{k-1}^2-}, \eqref{sum_{k=1}^n X_{k-2}^2-} and \eqref{sum_{k=1}^n X_{k-1} X_{k-2}-}
 yield \eqref{help_deterministic_scale1}.
\proofend

\appendix

\vspace*{5mm}

\noindent{\bf\Large Appendices}

\section{Stable convergence}
\label{section_stable}

We recall the notions of stable and mixing convergence.

\begin{Def}\label{Def_HL_stable_conv}
Let \ $(\Omega,\cF,\PP)$ \ be a probability space and \ $\cG\subset \cF$ \ be a sub-$\sigma$-field.
Let \ $(\bX_n)_{n\in\NN}$ \ and \ $\bX$ \ be \ $\RR^d$-valued random variables defined on \ $(\Omega,\cF,\PP)$,
 \ where \ $d\in\NN$.

\noindent (i) We say that \ $\bX_n$ \ converges \ $\cG$-stably to \ $\bX$ \ as \ $n\to\infty$, \ if the conditional distribution
 \ $\PP^{\bX_n\mid \cG}$ \ of \ $\bX_n$ \ given \ $\cG$ \ converges weakly to the conditional distribution
 \ $\PP^{\bX\mid \cG}$ \ of \ $\bX$ \ given \ $\cG$ \ as \ $n\to\infty$ \ in the sense of weak convergence of Markov kernels.
It equivalently means that
 \[
   \lim_{n\to\infty} \EE_\PP(\xi \EE_\PP(h(\bX_n) \mid \cG ) )
      = \EE_\PP( \xi \EE_\PP(h(\bX) \mid \cG ) )
 \]
 for all random variables \ $\xi:\Omega\to\RR$ \ with \ $\EE_\PP(\vert \xi\vert)<\infty$ \ and for all bounded and continuous functions
 \ $h:\RR^d\to\RR$.

\noindent (ii) We say that \ $\bX_n$ \ converges \ $\cG$-mixing to \ $\bX$ \ as \ $n\to\infty$, \
 if \ $\bX_n$ \ converges \ $\cG$-stably to \ $\bX$ \ as \ $n\to\infty$, \ and \ $\PP^{\bX\mid \cG} = \PP^\bX$ \ $\PP$-almost surely, where
 \ $\PP^\bX$ \ denotes the distribution of \ $\bX$ \ on \ $(\RR^d,\cB(\RR^d))$ \ under \ $\PP$.
\ Equivalently, we can say that \ $\bX_n$ \ converges \ $\cG$-mixing to \ $\bX$ \ as \ $n\to\infty$,
 \ if \ $\bX_n$ \ converges \ $\cG$-stably to \ $\bX$ \ as \ $n\to\infty$, \ and \ $\sigma(\bX)$ \ and \ $\cG$ \ are independent,
 which equivalently means that
 \[
   \lim_{n\to\infty} \EE_\PP(\xi \EE_\PP(h(\bX_n) \mid \cG ) )
      = \EE_\PP(\xi)\EE_\PP(h(\bX))
 \]
 for all random variables \ $\xi:\Omega\to\RR$ \ with \ $\EE_\PP(\vert \xi\vert)<\infty$ \ and for all bounded and continuous functions
 \ $h:\RR^d\to\RR$.
\end{Def}

In Definition \ref{Def_HL_stable_conv}, \ $\PP^{\bX_n\mid \cG}$, \ $n\in\NN$, \ and \ $\PP^{\bX\mid \cG}$ \ are
 the \ $\PP$-almost surely unique \ $\cG$-measurable Markov kernels from \ $(\Omega,\cF)$ \ to \ $(\RR^d,\cB(\RR^d))$ \ such that for each \ $n\in\NN$,
 \[
   \int_G \PP^{\bX_n\mid \cG}(\omega,B)\,\PP(\dd \omega)
     = \PP(\bX_n^{-1}(B) \cap G)
     \qquad \text{for every \ $G\in \cG$, \ $B\in\cB(\RR^d)$.}
 \]
 and
 \[
   \int_G \PP^{\bX\mid \cG}(\omega,B)\,\PP(\dd \omega)
     = \PP(\bX^{-1}(B) \cap G)
     \qquad \text{for every \ $G\in \cG$, \ $B\in\cB(\RR^d)$,}
 \]
 respectively.
For the notion of weak convergence of Markov kernels towards a Markov kernel,
 see H\"ausler and Luschgy \cite[Definition 2.2]{HauLus}.
For more details on stable convergence, see H\"ausler and Luschgy \cite[Chapter 3 and Appendix A]{HauLus}.

\section{A multidimensional stable limit theorem}
\label{section_MSLTES}

Recall that \ $\log^+(x):= \log(x)\bbone_{\{x\geq 1\}} + 0\cdot \bbone_{\{ 0\leq x < 1\}}$ \ for \ $x\in\RR_+$,
 \ and for an event \ $A$ \ with \ $\PP(A)>0$, \ $\PP_A$ \ denotes the conditional probability measure given \ $A$.
\ For an \ $\RR^d$-valued stochastic process \ $(\bU_n)_{n\in\ZZ_+}$, \ the increments
 \ $\Delta \bU_n$, \ $n \in \ZZ_+$, \ are defined by \ $\Delta \bU_0 := \bzero$ \ and
 \ $\Delta \bU_n := \bU_n - \bU_{n-1}$ \ for \ $n \in \NN$.

We recall a multidimensional analogue of Theorem 8.2 in H\"ausler and Luschgy \cite{HauLus} which was proved in Barczy and Pap \cite[Theorem 1.4]{BarPap1}.

\begin{Thm}[Barczy and Pap {\cite[Theorem 1.4]{BarPap1}}]\label{MSLTES}
Let \ $(\bU_n)_{n\in\ZZ_+}$ \ and \ $(\bB_n)_{n\in\ZZ_+}$ \ be \ $\RR^d$-valued and \ $\RR^{d\times d}$-valued stochastic processes, respectively,
 defined on a probability space \ $(\Omega,\cF,\PP)$ \ and
 adapted to a filtration \ $(\cF_n)_{n\in\ZZ_+}$.
\ Suppose that \ $\bB_n$ \ is invertible for
 sufficiently large \ $n \in \NN$.
\ Let \ $(\bQ_n)_{n\in\NN}$ \ be a sequence in \ $\RR^{d\times d}$ \ such that
 \ $\bQ_n \to \bzero$ \ as \ $n \to \infty$ \ and \ $\bQ_n$ \ is invertible for
 sufficiently large \ $n \in \NN$.
\ Let \ $G \in \cF_\infty := \sigma(\bigcup_{k=0}^\infty \cF_k)$ \ with \ $\PP(G) > 0$.
\ Assume that the following conditions are satisfied:
 \begin{enumerate}
  \item[\textup{(i)}]
   there exists an \ $\RR^{d\times d}$-valued, \ $\cF_\infty$-measurable random matrix \ $\Beta:\Omega\to\RR^{d\times d}$ \ such that
    \ $\PP(G \cap \{\exists\,\Beta^{-1}\}) > 0$ \ and
    \[
      \bQ_n \bB_n^{-1} \stochG \Beta \qquad \text{as \ $n \to \infty$,}
    \]
  \item[\textup{(ii)}]
   $(\bQ_n \bU_n)_{n\in\NN}$ \ is stochastically bounded in \ $\PP_{G\cap\{\exists\,\Beta^{-1}\}}$-probability, i.e.,
   \[
     \lim_{K\to\infty} \sup_{n\in\NN} \PP_{G\cap\{\exists\,\Beta^{-1}\}}(\|\bQ_n \bU_n\| > K) = 0.
   \]
  \item[\textup{(iii)}]
   there exists an invertible matrix \ $\bP \in \RR^{d\times d}$ \ with \ $\varrho(\bP) < 1$
    \ such that
    \[
      \bB_n \bB_{n-r}^{-1} \stochG \bP^r \qquad
      \text{as \ $n \to \infty$ \ for every \ $r \in \NN$,}
    \]
  \item[\textup{(iv)}]
   there exists a probability measure \ $\mu$ \ on \ $(\RR^d, \cB(\RR^d))$ \ with
    \ $\int_{\RR^d} \log^+(\|\bx\|) \, \mu(\dd\bx) < \infty$
    \ such that
    \[
      \EE_\PP\bigl(\ee^{\ii\langle\btheta,\bB_n\Delta\bU_n\rangle}
                     \mid \cF_{n-1}\bigr)
      \stochGeta \int_{\RR^d} \ee^{\ii\langle\btheta,\bx\rangle} \, \mu(\dd\bx)
      \qquad \text{as \ $n \to \infty$}
    \]
    for every \ $\btheta \in \RR^d$.
 \end{enumerate}
Then
 \begin{equation}\label{conv_BU}
  \bB_n \bU_n \to \sum_{j=0}^\infty \bP^j \bZ_j \qquad
  \text{$\cF_\infty$-mixing under \ $\PP_{G \cap \{\exists\,\Beta^{-1}\}}$ \ as \ $n \to \infty$,}
 \end{equation}
 and
 \begin{equation}\label{conv_QU}
   \bQ_n \bU_n \to \Beta \sum_{j=0}^\infty \bP^j \bZ_j \qquad
   \text{$\cF_\infty$-stably under \ $\PP_{G \cap \{\exists\,\Beta^{-1}\}}$ \ as \ $n \to \infty$,}
 \end{equation}
 where \ $(\bZ_j)_{j\in\ZZ_+}$ \ denotes a sequence of \ $\PP$-independent and identically distributed \ $\RR^d$-valued random vectors
 being \ $\PP$-independent of \ $\cF_\infty$ \ with \ $\PP(\bZ_0 \in B) = \mu(B)$ \ for all \ $B \in \cB(\RR^d)$.
\end{Thm}

The series \ $\sum_{j=0}^\infty \bP^j \bZ_j$ \ in \eqref{conv_BU} and in \eqref{conv_QU} is absolutely convergent \ $\PP$-almost surely,
 since, by condition (iv) of Theorem \ref{MSLTES}, \ $\EE_\PP(\log^+(\|\bZ_0\|))<\infty$ \ and one can apply Lemma 1.3 in Barczy and Pap \cite{BarPap1}.
Further, the random variable \ $\Beta$ \ and the sequence \ $(\bZ_j)_{j\in\ZZ_+}$ \ are \ $\PP$-independent in Theorem \ref{MSLTES},
 since \ $\Beta$ \ is \ $\cF_\infty$-measurable and the sequence \ $(\bZ_j)_{j\in\ZZ_+}$ \ is \ $\PP$-independent of \ $\cF_\infty$.

\section{Lenglart's inequality}
\label{section_Lenglart}

The following form of the Lenglart's inequality can be found, e.g., in H\"ausler and Luschgy \cite[Theorem A.8]{HauLus}.

\begin{Thm}\label{Lenglart_X}
Let \ $(\xi_n)_{n\in \ZZ_+}$ \ be a nonnegative submartingale with respect to a filtration \ $(\cF_n)_{n\in\ZZ_+}$ \ and with compensator
 \ $A_n := \sum_{k=1}^n \EE_\PP(\xi_k - \xi_{k-1} \mid \cF_{k-1})$, \ $n \in \NN$, \ $A_0:=0$.
\ Then for each \ $a, b \in \RR_{++}$ \ and \ $n\in\NN$,
 \[
   \PP\Big(\max_{k\in\{0,1,\ldots,n\}} \xi_k \geq a\Big)
   \leq \frac{b}{a}
        + \PP(\xi_0 + A_n > b) .
 \]
\end{Thm}

Applying Theorem \ref{Lenglart_X} for the square of a square integrable martingale, we obtain the following corollary.

\begin{Cor}\label{Lenglart_Y}
Let \ $(\eta_n)_{n\in \ZZ_+}$ \ be a square integrable martingale with respect to a filtration \ $(\cF_n)_{n\in\ZZ_+}$ \ and
 with quadratic characteristic process \ $\langle \eta\rangle_n := \sum_{k=1}^n \EE_\PP((\eta_k - \eta_{k-1})^2 \mid \cF_{k-1})$, \ $n \in \NN$,
 \ $\langle \eta\rangle_0:=0$.
\ Then for each \ $a, b \in \RR_{++}$ \ and \ $n\in\NN$,
 \[
   \PP\Big(\max_{k\in\{0,1,\ldots,n\}} \eta_k^2 \geq a\Big)
   \leq \frac{b}{a}
        + \PP(\eta_0^2 + \langle \eta\rangle_n > b) .
 \]
\end{Cor}
Note that under the conditions of Corollary \ref{Lenglart_Y}, the quadratic characteristic process \ $(\langle \eta\rangle_n)_{n\in\ZZ_+}$ \
 of \ $(\eta_n)_{n\in\ZZ_+}$ \ coincides with the compensator of \ $(\eta_n^2)_{n\in\ZZ_+}$.

\section*{Acknowledgements}
We are grateful to Michael Monsour for sending us his paper \cite{MikMon1} and the paper of Venkataraman \cite{Ven3}
 about the limiting distributions for the LSE of AR parameters of AR(2) processes.
This paper was finished after the sudden death of our longtime co-author, mentor
and dear friend Gyula Pap who passed away in October 2019. We dedicate this paper to him.
Also, we would like to thank the referee for her/his comments that helped us improve the paper.

\section*{Conflicts of interest}
On behalf of all authors, the corresponding author states that there is no conflict of interest.

\end{document}